\documentclass[12pt,twoside]{amsart}
\usepackage{mathrsfs}
\usepackage{}
\usepackage{txfonts}
\usepackage{url}
\usepackage{amsthm}
\usepackage{latexsym}
\usepackage{mathabx}
\usepackage[all]{xy}\SelectTips{eu}{}
\usepackage{amssymb}
\usepackage{amscd,mathrsfs,graphicx,mathrsfs}

\date{}
\allowdisplaybreaks[4] \footskip=15pt

\renewcommand{\uppercasenonmath}[1]{}

\numberwithin{equation}{section} \theoremstyle{plain}
\newtheorem*{thm*}{Main Theorem}
\newtheorem{thm}{Theorem}[section]
\newtheorem{cor}[thm]{Corollary}
\newtheorem*{cor*}{Corollary}
\newtheorem{lem}[thm]{Lemma}
\newtheorem*{lem*}{Lemma}

\newtheorem*{fact*}{Fact}

\newtheorem*{nota*}{Notation}
\newtheorem{prop}[thm]{Proposition}
\newtheorem*{prop*}{Proposition}
\newtheorem{rem}[thm]{Remark}
\newtheorem*{rem*}{Remark}

\newtheorem*{observation*}{Observation}
\newtheorem{exa}[thm]{Example}
\newtheorem*{exa*}{Example}
\newtheorem{df}[thm]{Definition}
\newtheorem*{df*}{Definition}

\newtheorem*{conj*}{Conjecture}

\newtheorem*{ques*}{Question}

\newtheorem*{ack*}{ACKNOWLEDGEMENTS}

\newcommand{\Z}{\mbox{\rm Z}}

\newcommand{\C}{\mbox{\rm C}}
\newcommand{\D}{\mbox{\rm D}}
\newcommand{\ac}{\mbox{\rm ac}}

\newcommand{\K}{\mbox{\rm K}}
\newcommand{\Ch}{\mbox{\rm Ch}}

\newcommand{\Hom}{\mbox{\rm Hom}}

\newcommand{\Ext}{\mbox{\rm Ext}}
\newcommand{\Tor}{\mbox{\rm Tor}}
\newcommand{\h}{\mbox{\rm H}}

\newcommand{\dg}{\mbox{\rm dg}}

\newcommand{\id}{\mbox{\rm id}}

\newcommand{\pd}{\mbox{\rm pd}}
\newcommand{\fd}{\mbox{\rm fd}}

\newcommand{\Gpd}{\mbox{\rm Gpd}}
\newcommand{\Gid}{\mbox{\rm Gid}}

\newcommand{\Gfd}{\mbox{\rm Gfd}}

\newcommand{\Mod}{\mbox{\rm Mod}}

\newcommand{\p}{\mbox{\rm Proj}}

\newcommand{\inj}{\mbox{\rm Inj}}

\newcommand{\im}{\mbox{\rm Im}}

\newcommand{\Ker}{\mbox{\rm Ker}}
\newcommand{\Coker}{\mbox{\rm Coker}}
\newcommand{\cone}{\mbox{\rm Cone}}

\newcommand{\Cone}{\mbox{\rm Cone}}

\newcommand{\Ggldim}{\mbox{\rm G-gldim}}
\newcommand{\Gwgldim}{\mbox{\rm G-wgldim}}
\newcommand{\gldim}{\mbox{\rm gldim}}
\newcommand{\wgldim}{\mbox{\rm wgldim}}

\newcommand{\GProj}{\mbox{\rm GProj}}
\newcommand{\GInj}{\mbox{\rm GInj}}
\newcommand{\Proj}{\mbox{\rm Proj}}
\newcommand{\Inj}{\mbox{\rm Inj}}

\newcommand{\Glgl}{\mbox{\rm Ggldim}}

\newcommand{\Gwgl}{\mbox{\rm Gwgldim}}

\newcommand{\sfli}{\mbox{\rm sfli}}
\newcommand{\FFD}{\mbox{\rm FFD}}

\newcommand{\Gfcd}{\mbox{\rm Gfcd}}

\begin{document}
\begin{center}
{\Large \bf Homotopy equivalences and Grothendieck duality over rings with finite Gorenstein weak global dimension }
\footnotetext{* Corresponding author.

E-mail:~wangjunpeng1218@163.com; sestrada@um.es}

\vspace{0.5cm} {\small Junpeng Wang$^{1}$ and Sergio Estrada$^{2,*}$ \\\

1.  \emph{Department of Mathematics}, \emph{Northwest Normal University}, \\ \emph{Lanzhou 730070, People's Republic of China} \\
}
2. \emph{Department of Mathematics}, \emph{Universidad de Murcia}, \\ \emph{Murcia 30100, Spain}

\end{center}

\vskip.5cm

\noindent{\bf Abstract:}
{Let $R$ be a ring with Gwgldim$(R)<\infty$. We obtain a triangle-equivalence
$\K(R\text{-}\GProj)\simeq \K(R\text{-}\GInj)$
which restricts to a triangle-equivalence $\K(R\text{-}\Proj)$ $\simeq \K(R\text{-}\Inj)$. This class of rings includes, among others, (left) Gorenstein rings, Ding-Chen rings and the more general Gorenstein $n$-coherent rings ($n\in \mathbb{N}\cup \{\infty\}, n\geq 2$).
As application, we establish some triangle-equivalences of Grothendieck duality over Ding-Chen rings and Gorenstein $n$-coherent rings.
}

\vskip.2cm

 \vskip.2cm

 \vspace{10pt} \vspace{10pt}

\section{\bf Introduction}

Gorenstein rings were introduced by Iwanaga in \cite{Iwa1979}.
They play an important role in (Gorenstein) homological algebra, representation theory of groups and algebras, etc.
As a generalization of Gorenstein rings, \emph{Ding-Chen rings} were initially introduced by Ding and
Chen in \cite{DC1996} under the name ``$m$-FC rings'' and renamed by Gillespie in \cite{Gil2010}.
To be more general, Wang, Liu, and Yang introduced \emph{ Gorenstein} $n$-\emph{coherent rings} in \cite{WLY2018}.
It is defined as a two-sided $n$-coherent ring $R$ such that all FP$_n$-injective left and right $R$-modules have finite flat dimension at most $m$ for some integer $m\geq 0$ (see Definition \ref{9-1}).
One motivation for introducing such rings is to give a negative answer to a question posed by Gillespie in \cite{Gil2017}.

Another generalization of Gorenstein rings is the so-called \emph{ left Gorenstein ring}, introduced by Beligiannis \cite{Bel2000}. It is defined as a ring $R$ such that all projective left $R$-modules have finite injective dimension and all injective left $R$-modules have finite projective dimension. This ring is just a ring whose left Gorenstein global dimension, Ggldim$(-)$, is finite (see, for example, \cite[Theorem 4.1]{Emm2012}). In addition, by \cite[Corollary 6.11]{Bel2000},  a ring is Gorenstein if and only if it is two-sided noetherian and left Gorenstein.  One important feature about left Gorenstein rings is that they provide canonical triangle equivalences between well-known homotopy categories. As usual, denote by $\K(R\text{-}\GProj)$ (resp. $\K(R\text{-}\GInj)$) the homotopy category of Gorenstein projective (resp. Gorenstein injective) left $R$-modules.
Its homotopy subcategory consisting of projective (resp. injective) left $R$-modules is denoted by $\K(R\text{-}\Proj)$ (resp. $\K(R\text{-}\Inj)$).
We also denote by $\K_{\mathrm{ac}}(R\text{-}\Proj)$ (resp. $\K_{\mathrm{ac}}(R\text{-}\Inj)$) the homotopy subcategory of $\K(R\text{-}\Proj)$ (resp. $\K(R\text{-}\Inj)$) consisting of exact complexes of projective (resp. injective) left $R$-modules.
The study of compactness of $\K(R\text{-}\Proj)$ and $\K(R\text{-}\Inj)$ and triangle equivalence between them was initiated in 2005 by Krause \cite{Kra2005} and J\o rgensen~\cite{Jor2005}, and has been extensively investigated during the last years by Iyengar, Krause, Neeman, Chen, \v{S}t'ov\'{\i}\v{c}ek, Asadollahi, Hafezi, Salarian, Positselski and Gillespie \cite{IK2006, Nee2008, Nee2014, Chen2010, Sto2014, AHS2014, Gil2017, Pos2017}. Now the relation between such homotopy categories and left Gorenstein rings is made explicit in the following result due to Chen

\begin{thm} \label{02} {\rm (=\cite[Theorem B]{Chen2010}).}
Let $R$ be a left Gorenstein ring.
 Then there is a triangle-equivalence
$\K(R\text{-}\GProj)\simeq \K(R\text{-}\GInj)$
which restricts to a triangle-equivalence $\K(R\text{-}\Proj)\simeq \K(R\text{-}\Inj).$
 \end{thm}

Another classical Gorensteing homological global invariant is the Gorenstein weak global dimension of a ring $R$, Gwgldim$(R)$. It is a refinement of the usual weak global dimension of $R$,  and it is closely related to Ggldim$(R)$.
Recently, we have used different methods to prove that Gwgldim$(R)\leq $ Ggldim$(R)$ in \cite[Theorem 3.7]{WYZ2023} and in \cite[Theorem 5.13]{C-W2023}. In other words, any left Gorenstein ring admits finite Gorenstein weak global dimension.
On the other hand, it is shown in Remark \ref{06-3-1} that a ring $R$ is of Gwgldim$(R)< \infty$ if and only if all injective left and right $R$-modules have finite flat dimension.
Thus, it is easy to see that any Gorenstein $n$-coherent ring $R$ (in particular, any Gorenstein and Ding-Chen ring) admits  Gwgldim$(R)< \infty$.

Our main result substantially improves Theorem \ref{02} by relaxing the assumption of the finiteness of the Gorenstein global dimension to the weak one. Notice that there are examples of rings with infinite Gorenstein global dimension but finite Gorenstein weak global dimension (see Example \ref{9-6}).

\begin{thm} \label{03} {\rm (=Theorem \ref{8}).}
Let $R$ be a ring with  $\mathrm{Gwgldim}(R)< \infty$.
 Then there is a triangle-equivalence
$\K(R\text{-}\GProj)\simeq \K(R\text{-}\GInj)$
which restricts to a triangle-equivalence $\K(R\text{-}\Proj)\simeq \K(R\text{-}\Inj).$
 \end{thm}

We also refer to \cite{Chen2010} to note the following: it is not clear, just using Theorem \ref{02} and its proof there, that over a \emph{left} Gorenstein ring $R$ whether or not the following hold:

(1) the triangle equivalence $\K(R\text{-}\GProj)\simeq \K(R\text{-}\GInj)$
can restrict to a triangle-equivalence $$\K_{\mathrm{ac}}(R\text{-}\Proj)\simeq \K_{\mathrm{ac}}(R\text{-}\Inj);\quad\text{and}$$

(2) the triangle equivalences involving the corresponding homotopy categories of \emph{right} $R$-modules hold true as well.

With our proof of Theorem \ref{03}, we give an explicit construction of the equivalent functor between $\K(R\text{-}\GProj)$ and $\K(R\text{-}\GInj)$, which guarantees that the answer to (1) is affirmative (see Corollary \ref{9-0-0}).
On the other hand, it is known from \cite[Corollary 2.5]{CET2021} that the notion of Gorenstein weak global dimension of rings is left-right symmetric. Hence the answer to (2) is also affirmative (see Corollary \ref{9-0-0}).

The application of Theorem \ref{03} is twofold. One is to obtain a triangle equivalence which relates the subcategory $\mathcal{A}$ and the the subcategory of flat modules whenever $R$ is a flat-typed $(\mathcal{L}, \mathcal{A})$-Gorenstein ring with respect to a bi-complete duality pair $(\mathcal{L}, \mathcal{A})$ (see Theorem \ref{7-5} for the result and see Definitions 2.1-2.3 for the corresponding notions).
Moreover, we show that the equivalent triangulated categories
are compactly generated under some mild conditions (see Theorem \ref{7-7}).

The other is to obtain some triangle equivalences involving Grothendieck duality over Ding-Chen rings and more general Gorenstein $n$-coherent rings.
Let us give a simple background of Grothendieck duality and describe our work in more detail.

Grothendieck duality is a classical subject which can go back 1958s. Roughly speaking, it is a statement concerning the existence of a right adjoint to
the ``direct image with compact support'' functor between derived categories of sheaves or modules.
For a ring $R$, over the years many people investigated Grothendieck duality for derived categories of $R$-modules by providing the following insights:

\begin{enumerate}
\item[(\textbf{GD1})] $R$ has a dualizing complex.
\item[(\textbf{GD2})] There is a triangle-equivalence $\mathbf{D}^{b}(R^{op}\text{-}\mathrm{mod})^{op}\simeq \mathbf{D}^{b}(R\text{-}\mathrm{mod}),$
where the right category is the bounded derived category of finitely presented left $R$-modules and the left one is the opposite category of the bounded derived category of finitely  presented right $R$-modules;
  \item[(\textbf{GD3})] $\K(R\text{-}\p)$ and $\K(R\text{-}\inj)$ are triangulated equivalent.
 \end{enumerate}
According to \cite[Introduction]{Nee2008}, we know that,  over a left noetherian and right coherent ring $R$ such that all flat left $R$-modules have finite projective dimension, there are implications
$$(\textbf{GD1})\Longrightarrow (\textbf{GD3}) \Longrightarrow(\textbf{GD2}).$$ In Remark \ref{10-12} we note that the assumptions on $R$ to fulfill these implications can be relaxed to coherent rings $R$ such that every FP-injective left $R$-module has finite injective dimension.
Also, for a left noetherian and right coherent ring $R$, Theorem \ref{02} and its corollary ~\cite[Corollary C]{Chen2010} show that the above implications remain true if we replace $(\textbf{GD1})$ with ``$R$ is a left Gorenstein ring''. We note in Remark \ref{10-10} that the left noetherian assumption for $R$ can be relaxed to \emph{left coherent}.

Motivated by these results, if we apply Theorem \ref{03} to Ding-Chen rings, we obtain that, over a left and right coherent ring $R$,  the above implications still hold true if we replace $(\textbf{GD1})$ with ``$R$ is a Ding-Chen ring'' (see Corollary \ref{10}).




Let $R$ be a left noetherian ring. Denote by $\K^b (R$-proj) the bounded homotopy category of finitely presented and projective left $R$-modules and
let $\mathbf{D}^b(R$-mod) be as above. The Verdier quotient triangulated category
$$\mathbf{D}_{sg}(R):= \mathbf{D}^b(R\text{-}\mathrm{mod})/\K^b (R\text{-}\mathrm{proj})$$
is known as the \emph{singularity category} of left $R$-modules in the literature.
This Verdier quotient was introduced primarily by Buchweitz in \cite{Buc1986} under the name ``stable derived category'' and renamed by Orlov in \cite{Orl2004}.
On the one hand, by Buchweitz's theorem (see \cite[Theorem 4.1]{Buc1986}) singularity categories are related with Gorenstein projective modules.
On the other hand,  using Krause's recollement,
Krause \cite[Corollary 5.4]{Kra2005} proved that the homotopy category $\K_{\mathrm{ac}}(R\text{-}\Inj)$ is compactly generated
with a triangle equivalence
$$\K_{\mathrm{ac}}(R\text{-}\Inj)^{c}\simeq \mathbf{D}_{sg}(R),$$
where  $\mathcal{T}^{c}$ denotes the triangulated subcategory of compact objects of a triangulated category $\mathcal{T}$.




Motivated by \cite[Corollary 5.4]{Kra2005} and applying Theorem \ref{03} to Ding-Chen rings and left Gorenstein rings,
we obtain the following. 

\begin{cor} \label{04} {\rm (=Corollary \ref{12-1}).}
Let $R$ be a left and right coherent ring. Consider the following conditions
\begin{enumerate}
\item[(\textbf{GD1$'$})] $R$ is Ding-Chen or left Gorenstein;
\item[(\textbf{GD2$'$})] There is a triangle-equivalence $\mathbf{D}_{sg}(R^{op})^{op}\simeq \mathbf{D}_{sg}(R),$
where the left category is the opposite category of the singularity category of right $R$-modules;
  \item[(\textbf{GD3$'$})] $\K_{\mathrm{ac}}(R\text{-}\p)$ and $\K_{\mathrm{ac}}(R\text{-}\inj)$ are triangulated equivalent.
 \end{enumerate}
Then there are implications
$(\textbf{GD1$'$})\Longrightarrow (\textbf{GD3$'$}) \Longrightarrow(\textbf{GD2$'$}).$
\end{cor}

Finally, the above two results for Ding-Chen rings have both a version for Gorenstein $n$-coherent rings (see Corollary \ref{11} and Theorem \ref{12}).

\section{\bf Preliminaries}

Throughout this article, all rings $R$ are assumed to be associative with identity and all modules
are unitary.

In this section, we mainly recall some necessary notions and facts, which will be used in the paper.
Let $R$ be a ring. The article will involve left and right $R$-modules; complexes of left and right $R$-modules;
homotopy categories of left and right $R$-modules and derived categories of left and right $R$-modules.

We denote by $R$-Mod the category of all left $R$-modules. Its subcategory consisting of all projective (resp. injective, flat,  FP-injective and cotorsion modules)
modules will be denoted by $_R\mathcal{P}$ (resp. $_R\mathcal{I}$, $_R\mathcal{F}$, $_R\mathcal{FI}$ and $_R \mathcal{C}ot$). For the corresponding (sub)category of right $R$ modules, we use the notation
$R^{op}$-Mod or Mod-$R$ to stand for the category of all right $R$-modules, and use $_{R^{op}}\mathcal{P}$ or $\mathcal{P}_R$ to denote the subcategory consisting of all projective right $R$-modules, etc.,
where $R^{op}$ means the opposite ring of $R$. In addition, the notation $\mathcal{P}$ (resp. $\mathcal{I}$ and $\mathcal{F}$) will denote the subcategory
 consisting of all projective (resp. injective and flat) left or right (not necessarily left and right) $R$-modules. Denote by $\pd(_RM)$ (resp. $\id(_RM)$, $\fd(_RM)$) the projective (resp. injective, flat) dimension of a left $R$-module $_RM$;
by $\gldim(R)$ (resp. $\gldim(R^{op})$, $\wgldim(R)$) the left global (right global, weak global) dimension of $R$.

We denote by $\Ch(R\text{-}\Mod)$ (resp. $\K(R\text{-}\Mod)$) the category of complexes of left $R$-modules (resp. the homotopy category of left $R$-modules).
Their subcategory consisting of all projective (resp. injective, Gorenstein projective and Gorenstein injective) complexes of left $R$-modules
will be denoted by $\Ch(R\text{-}\Proj)$, $\Ch(R\text{-}\Inj)$, $\Ch(R\text{-}\GProj)$ and $\Ch(R\text{-}\GInj)$
as well as $\K(R\text{-}\Proj)$, $\K(R\text{-}\Inj)$, $\K(R\text{-}\GProj)$ and $\K(R\text{-}\GInj)$, respectively.
For the corresponding (sub)category of right $R$-modules, we use the notations $\Ch(\Mod\text{-}R)$, $\K(\Mod\text{-}R)$, $\Ch(\Proj\text{-}R)$ and $\K(\Proj\text{-}R)$ and so on.
Furthermore, for a subcategory $_R\mathcal{X}^{\bullet}$ of $\Ch(R\text{-}\Mod)$, we denote by $\K(_R\mathcal{X}^{\bullet})$ the corresponding homotopy subcategory of $\K(R\text{-}\Mod)$.
Given a subcategory $_R\mathcal{X}$ of $R\text{-}\Mod$, we denote by $\Ch(_R\mathcal{X})$ the subcategory of $\Ch(R\text{-}\Mod)$ consisting of all complexes with all components in $_R\mathcal{X}$,
and simply denote $\K(\Ch(_R\mathcal{X}))$ by $\K(_R\mathcal{X})$. When considering exact (acyclic) complexes, we use the notations $\Ch_{\mathrm{ac}}(-)$ and $\K_{\mathrm{ac}}(-)$.

We denote by $\mathbf{D}(R\text{-}\Mod)$ the derived category of all left $R$-modules.
Given a subcategory $_R\mathcal{X}$ of $R\text{-}\Mod$, we denote by $\mathbf{D}_{\mathcal{X}}(R\text{-}\Mod)$ the relative derived category with respect to $_R\mathcal{X}$ (for the details, see Definition \ref{10-9-2}).
In particular, we denote by $\mathbf{D}_{\mathrm{GP}}(R\text{-}\Mod)$ (resp. $\mathbf{D}_{\mathrm{GI}}(R\text{-}\Mod)$) the relative derived category with respect to Gorenstein projective (resp. Gortenstein injective)
left $R$-modules. For the corresponding (sub)category of right $R$-modules, the notations $\mathbf{D}(R^{op}\text{-}\Mod)$, $\mathbf{D}_{\mathrm{GP}}(R^{op}\text{-}\Mod)$ and  $\mathbf{D}_{\mathrm{GI}}(R^{op}\text{-}\Mod)$ are used.

\subsection{\bf  Cotorsion pairs.}

 Let $\mathcal{C}$ be an abelian category and $\mathcal{X},\mathcal{Y}, \mathcal{Z}$ be subcategories of $\mathcal{C}$.
A pair $(\mathcal{X},\mathcal{Y})$ is called a \emph{cotorsion pair} if
$\mathcal{X}^\perp=\mathcal{Y}$ and $\mathcal{X}=$
$^\perp\mathcal{Y}$. Here $\mathcal{X}^\perp=$ $\{C\in \mathcal{C}~|~\mathrm{Ext}^1_{\mathcal{C}}(X,C)=0,
 \forall X\in\mathcal{X}\}$, and  $^\perp\mathcal{Y}$ can be defined dually.
A cotorsion pair $(\mathcal{X},\mathcal{Y})$ is said to be \emph{hereditary} if Ext$^m_{\mathcal{C}}(X,Y)=0$ for all $X\in \mathcal{X}$, $Y\in \mathcal{Y}$ and $m\geq 1$.
A cotorsion pair $(\mathcal{X},\mathcal{Y})$ is said to be $cogenerated~by~a~set$ if there is a set $\mathcal{S}\subseteq \mathcal{X}$ of objects in $\mathcal{C}$ such that $\mathcal{Y}=\mathcal{S}^\perp$.
A cotorsion pair $(\mathcal{X},\mathcal{Y})$ is called \emph{complete} if
for any object $C$ of $\mathcal{C}$, there are short exact
sequences $0\rightarrow Y\rightarrow X\rightarrow C\rightarrow 0$ and $0\rightarrow C\rightarrow Y'\rightarrow X'\rightarrow 0$ in $\mathcal{C}$ with
  $X,X'\in\mathcal{X}$ and $Y,Y'\in\mathcal{Y}$.
A triple $(\mathcal{X},\mathcal{Z},\mathcal{Y})$ is called a \emph{cotorsion triple} if the two pairs $(\mathcal{X},\mathcal{Z})$ and $(\mathcal{Z},\mathcal{Y})$ are cotorsion pairs.
We say that a cotorsion triple $(\mathcal{X},\mathcal{Z},\mathcal{Y})$ is \emph{complete} (resp. \emph{hereditary})
if the two cotosion pairs $(\mathcal{X},\mathcal{Z})$ and $(\mathcal{Z},\mathcal{Y})$ are complete (resp. hereditary).

Given $C\in \mathcal{C}$ and  $\mathcal{X} \subseteq$ $\mathcal{C}$, a \emph{special $\mathcal{X}$-preenvelope} of $C$ is defined as a monic
morphism $\alpha : C \to X$ with $X \in \mathcal{X}$ and $\Coker\alpha \in {^{\perp}\mathcal{X}}$.
A subcategory $\mathcal{X}$ of $\mathcal{C}$ is called \emph{injectively coresolving} if it is closed under extensions and cokernels of monic morphisms, and contains all injective objects of $\mathcal{C}$.
Dually, we have the definitions of
\emph{special} $\mathcal{X}$-\emph{precover} and that $\mathcal{X}$ is  \emph{projectively resolving}.

We say that a subcategory $\mathcal{X}$ of $\mathcal{C}$ is \emph{thick} if $\mathcal{X}$ is closed under direct summands and such that if any two out of three of the terms in
a short exact sequence are in $\mathcal{X}$, then so is the third.

If $\mathcal{C}$ has enough projectives and enough injectives, then we say that a pair $(\mathcal{X},\mathcal{Y})$ of subcategories of $\mathcal{C}$ is a \emph{projective cotorsion pair} (resp. \emph{injective cotorsion pair}) if $(\mathcal{X},\mathcal{Y})$ forms a complete cotorsion pair such that $\mathcal{Y}$ (resp. $\mathcal{X}$) is thick and $\mathcal{X}\cap\mathcal{Y}$ coincides with the subcategory of all projective (resp. injective) objects. Equivalently, if $(\mathcal{X},\mathcal{Y})$ forms a hereditary complete cotorsion pair such that $\mathcal{X}\cap\mathcal{Y}$ coincides with the subcategory of all projective (resp. injective) objects (see \cite[Definition 3.4]{Gil2016} and \cite[Propositions 3.6 and 3.7]{Gil2016}).

\subsection{\bf Duality pairs and relative Gorenstein rings}

As usual, we write
$(_RM)^{+}=\Hom_\mathbb{Z}(_RM,\mathbb{Q}/\mathbb{Z})$ for a left $R$-module $_RM$.
The following definitions can be found in \cite[Definitions 3.1, 3.2 and 4.2]{WD2020}.

\begin{df} \label{06} {\rm Let $R$ be a ring.

(1) A \emph{duality pair} over $R$ is a pair $(_R\mathcal{L},\mathcal{A}_R)$, where $_R\mathcal{L}$ (resp. $\mathcal{A}_R$ ) is a subcategory consisting of left (resp. right)
$R$-modules, subject to the following two conditions:

(i) \,$_RL \in ~_R\mathcal{L}$ if and only if $(_RL)^{+} \in ~\mathcal{A}_R$, and

(ii) $\mathcal{A}_R$ is closed under direct summands and finite direct sums.

(2) A duality pair $(_R\mathcal{L},\mathcal{A}_R)$ over $R$ is called \emph{(co)product-closed} if the class
$_R\mathcal{L}$ is closed under arbitrary (co)products.

(3) A duality pair $(_R\mathcal{L},\mathcal{A}_R)$ over $R$ is called \emph{perfect} if $(_R\mathcal{L},\mathcal{A}_R)$  is coproduct-closed, and if $_R\mathcal{L}$ is
closed under extensions and contains the regular module $_RR$.

(4) A duality pair $(_R\mathcal{L},\mathcal{A}_R)$ over $R$ is called  \emph{symmetric} if $(_R\mathcal{L},\mathcal{A}_R)$ is a duality pair
over $R$ and $(_{R^{op}}\mathcal{A},\mathcal{L}_{R^{op}})$ is a duality pair over $R^{op}$.

(5) A duality pair $(_R\mathcal{L},\mathcal{A}_R)$ over $R$ is called  \emph{complete} if $(_R\mathcal{L},\mathcal{A}_R)$ is a symmetric and
perfect duality pair over $R.$

(6) A duality pair $(_R\mathcal{L},\mathcal{A}_R)$ over $R$ is called  \emph{bi-complete} if $(_R\mathcal{L},\mathcal{A}_R)$
is a complete duality pair over $R$ such that the pair $(^{\perp}\mathcal{A}_{R},\mathcal{A}_R)$ is a hereditary cotorsion pair in $R$-Mod which is
cogenerated by a set.
}
\end{df}

\begin{df} \label{07} {\rm
Let ($\mathcal{L}$, $\mathcal{A}$) be a pair of subcategories of left or right $R$-modules.
According to \cite[p.11]{WD2020},  we say that $(\mathcal{L},\mathcal{A})$ is a \emph{duality pair} over $R$ if $(_R\mathcal{L},\mathcal{A}_R)$ is a duality pair over $R$ and $(_{R^{op}}\mathcal{L}, \mathcal{A}_{R^{op}})$ is a duality pair over $R^{op}$. Similarly, we have the notions of that the duality pair $(\mathcal{L},\mathcal{A})$ over $R$ is \emph{(bi)-complete}.
}\end{df}

\begin{df} \label{2} {\rm
Let $R$ be a ring and ($\mathcal{L}$, $\mathcal{A}$) be a pair of subcategories of left or right $R$-modules.

(1) For a nonnegative integer $m$, $R$ is called an \emph{$m$-$(\mathcal{L},\mathcal{A})$}-\emph{Gorenstein ring} (resp. \emph{flat-typed $m$}-$(\mathcal{L},\mathcal{A})$-\emph{Gorenstein ring})
if the pair $(\mathcal{L}, \mathcal{A})$ is a bi-complete duality over $R$ such that all modules in $\mathcal{A}_{R^{op}}$ and $\mathcal{A}_{R}$ have finite $\mathcal{L}$-projective (resp. flat dimension) at most $m$ for some nonnegative integer $m$.

(2) $R$ is called an \emph{$(\mathcal{L},\mathcal{A})$}-\emph{Gorenstein ring} (resp. a \emph{flat-typed}~$(\mathcal{L},\mathcal{A})$-\emph{Gorenstein ring})
if there exists a nonnegative integer $m$ such that $R$ is an $m$-$(\mathcal{L},\mathcal{A})$-Gorenstein ring (resp. flat-typed $m$-$(\mathcal{L},\mathcal{A})$-Gorenstein ring).
}
\end{df}

\begin{rem} \label{2-0-0} {\rm
(1) The notion of duality pair was initially introduced by Holm and J{\o}rgensen in \cite{HJ2009} (for an arbitrary ring).

(2) The notion of complete (resp. symmetric) duality pairs was introduced by Gillespie in \cite[Definition 2.4]{Gil2019} for commutative rings.
The Definition \ref{06}(4)-(5) can be viewed as an arbitrary-ring version of \cite[Definition 2.4]{Gil2019}.

(3) Although Definition \ref{07} was not written as a ``Definition'' in \cite{WD2020} (just as a remark in \cite[p.11]{WD2020}),
it is important for the notion of (flat-typed)~$(\mathcal{L},\mathcal{A})$-Gorenstein rings.

(4) It is obvious that  any flat-typed ~$(\mathcal{L},\mathcal{A})$-Gorenstein  ring is always an ~$(\mathcal{L},\mathcal{A})$-Gorenstein ring.
As explained in \cite[Remark 4.4]{WD2020},  $(\mathcal{L},\mathcal{A})$-Gorenstein rings
are the same as~flat-typed ~$(\mathcal{L},\mathcal{A})$-Gorenstein rings in some cases, and there exist $(\mathcal{L},\mathcal{A})$-Gorenstein rings
are not a flat typed $(\mathcal{L},\mathcal{A})$-Gorenstein ring in other cases.
}
\end{rem}

\subsection{\bf Some basic notions of complexes.}

 A complex $_RX^{\bullet}$ in $\Ch(R\text{-}\Mod)$ is a sequence of left $R$-modules
$$_RX^{\bullet}=\xymatrix@C=0.5cm{
\cdots \ar[r] & _RX^{k-1} \ar[rr]^{d^{k-1}} && _RX^k \ar[rr]^{d^k} && _RX^{k+1} \ar[r]^{} & \cdots}$$
with $d^{k}d^{k-1}=0$ for all $k\in \mathbb{Z}$.
Given a complex $_RX^{\bullet}$,
the $k$th \emph{cycle} (resp., \emph{homology}) module of $_RX^{\bullet}$
is defined as Ker$d^k$ (resp., Ker$d^k/\im d^{k-1}$)
and denoted by $\Z^k(_RX^{\bullet})$ (resp., H$^{k}(_RX^{\bullet})$). We also let $\C^k(_RX^{\bullet})=$ Coker$d^{k-1}$.
Let
$$ \sup \h (_RX^{\bullet}) = \sup\{\ l\in\mathbb{Z} \ | \ \h^l(X) \neq0\}\quad\text{and} \quad \inf \h (_RX^{\bullet}) =\inf\{\ l\in\mathbb{Z}\  | \ \h^l(X)\neq0\}.$$
We follow the convention that $\sup \h (X) = -\infty$ and $\inf \h (X) = \infty$
if $X$ is \emph{exact} or \emph{acyclic} (i.e., $\h^k(_RX^{\bullet}) = 0$ for all $k\in \mathbb{Z}$).
The complex $_RX^{\bullet}$ is called \emph{homology bounded} (resp., \emph{bounded})
if H$^k(_RX^{\bullet}) = 0$ (resp., $_RX^{k}= 0$) for $|k| \gg 0$.

A \emph{morphism} $f^{\bullet}: ~_RX^{\bullet}\to ~_RY^{\bullet}$ of complexes is a family of morphisms
$f^{\bullet} = (f^{k}:~_RX^{k}\to ~_RY^{k})_{k\in \mathbb{Z}}$ of modules
satisfying $d_{Y^{\bullet}}^{k}f^{k} = f^{k+1}d^k_{X^{\bullet}}$ for each $k$.
A morphism $f^{\bullet}: ~_RX^{\bullet}\to ~_RY^{\bullet}$ of complexes
is called \emph{null-homotopic} if there is a family ~$(h^k:~_RX^{k}\to ~_RY^{k-1})_{k\in\mathbb{Z}}$ of morphisms of left $R$-modules
such that $$f^{k} = d_{Y^{\bullet}}^{k-1}h^{k} + h^{k+1}d_{X^{\bullet}}^{k}, \quad\forall k\in\mathbb{Z}.$$
A morphism $f^{\bullet}: ~_RX^{\bullet}\to ~_RY^{\bullet}$ of complexes
is said to be  a \emph{quasi-isomorphism} if the induced map $\h^k(f^{\bullet}):\h^k(_RX^{\bullet})\to \h^k(_RY^{\bullet})$ is an isomorphism for all $k\in \mathbb{Z}$.

Given a morphism~$f^{\bullet}: ~_RX^{\bullet}\rightarrow ~_RY^{\bullet}$,
the \emph{mapping cone} $\Cone(f^{\bullet})$ of $f^{\bullet}$ is defined as
a complex in $\Ch(R\text{-}\Mod)$ given by
$$\cone(f^{\bullet})^{k} = ~_RY^{k}\oplus ~_RX^{k+1} \quad \text{and} \quad
d_{\mathrm{Cone}(f^{\bullet})}^{k} = \left(
                                                 \begin{array}{cc}
                                                   d_{Y^{\bullet}}^k & f^{k+1} \\
                                                   0 & -d_{X^{\bullet}}^{k+1} \\
                                                 \end{array}
                                               \right), \quad \forall k\in\mathbb{Z}.$$
The \emph{homomorphism complex} of $_RX^{\bullet}$ and $_RY^{\bullet}$, denoted by $\mathcal{H}om(_RX^{\bullet},~_RY^{\bullet})$,
is defined as the complex
$$\cdots\longrightarrow \prod_{k\in \mathbb{Z}}\textrm{Hom}_{R}(_RX^{k}, ~_RY^{k+m}) \stackrel{d^{m}} \longrightarrow \prod_{k\in \mathbb{Z}}
\textrm{Hom}_{R}(_RX^{k}, ~_RY^{k+m+1})\longrightarrow\cdots$$
in $\Ch(\mathbb{Z}\text{-}\Mod)$, where $d^{m}$ is given by $d^{m}((f^{k})_{k\in \mathbb{Z}})=(d^{k+m}_{Y^{\bullet}}f^{k}-(-1)^{m}f^{k+1}d^{k}_{X^{\bullet}})_{k\in \mathbb{Z}}$, for all $(f^{k})_{k\in \mathbb{Z}} \in \mathcal{H}om(_RX^{\bullet}, ~_RY^{\bullet})^{m}.$

\subsection{\bf Gorenstein homological modules and dimensions}

Recall that a left $R$-module $_RM$ is
\emph{Gorenstein projective} \cite{EJ1995} if there exists  a totally acyclic complex $_RP^{\bullet}\in \Ch(R\text{-Proj})$
such that $_RM\cong \Z^{0}(_RP^{\bullet})$.
Here we say that a complex $_RP^{\bullet}\in \Ch(R\text{-Proj})$ is \emph{totally acyclic} if $_RP^{\bullet}$ is exact and $\Hom_R(_RP^{\bullet},~_RQ)$
remains exact for all projective left $R$-modules $_RQ$.
The notions of \emph{Gorenstein injective modules} and \emph{totally acyclic complexes of injective modules}
are defined dually.

A left $R$-module $_RM$ is said to be \emph{Gorenstein flat} \cite{EJT1993}
if there exists an \textbf{F}-totally acyclic complex $_RF^{\bullet}$ of flat left $R$-modules
such that $_RM\cong \Z^{0}(_RF^{\bullet})$.
Here we say that a complex $_RF^{\bullet}$ of flat left $R$-modules is \textbf{F}-\emph{totally acyclic} if $_RF^{\bullet}$ is exact and $I_R\otimes_R ~_RF^{\bullet}$
remains exact for all injective right $R$-modules $I_R$.

We denote by $_R\mathcal{GP}$ (resp., $_R\mathcal{GI}$ and $_R\mathcal{GF}$) the subcategory of $R\text{-}\Mod$ consisting of all
Gorenstein projective (resp., Gorenstein injective and Gorenstein flat) left $R$-modules. The corresponding subcategories of right $R$-modules and left or right $R$-modules are denoted by, $\mathcal{GP}_R$, $\mathcal{GI}_R$,  $\mathcal{GF}_R$ and $\mathcal{GP}$, $\mathcal{GI}$,  $\mathcal{GF}$, respectively.

Let $\mathcal{C}$ be an abelian category,
$\mathcal{X}$ its subcategory, and $C$ a object of $\mathcal{C}$.
Then the $\mathcal{X}$-\emph{projective dimension} of $C$, denoted by $\mathcal{X}\text{-}\pd(C)$,
is defined as follows:
$$ \mathcal{X}\text{-}\pd(C) = \inf \left \{m\in \mathbb{Z} \hspace{0.2cm}
\begin{array}{ |l} \xymatrix@C=6mm{0\ar[r]&X_m\ar[r]&\cdots\ar[r]&X_1\ar[r]&X_0\ar[r]&M\ar[r]&0} \text{is } \\
\text{ an exact sequence} \text{\ of objects in \ }\mathcal{C} \text{ with each } X_i\in\mathcal{X}.\end{array}\right\}
$$
If such an exact sequence does not exist, then we set $\mathcal{X}\text{-}\pd(C)=\infty$.
The $\mathcal{X}$-\emph{injective dimension} of $M$, $\mathcal{X}\text{-}\id(C)$, is defined dually.

In particular, the \emph{Gorenstein projective}
(resp., \emph{Gorenstein injective} and \emph{Gorenstein flat})
\emph{dimension} of a left $R$-module $_RM$ is defined as the $_R\mathcal{GP}$-projective (resp., $_R\mathcal{GI}$-injective and $_R\mathcal{GF}$-projective) dimension of $_RM$, and is denoted by
$\Gpd(_RM)$ (resp., $\Gid(_RM)$ and $\Gfd(_RM)$).

As a refinement of the usual global (resp., weak global) dimension of rings,
Gorenstein global (resp., Gorenstein weak global) dimension of rings is defined as follows:

\begin{df} \label{06-2}
{\rm For any ring $R$, its \emph{left Gorenstein} (resp. \emph{Gorenstein weak}) \emph{global dimension}, denoted by Ggldim$(R)$ (resp. Gwgldim$(R)$),  is defined via the following formula
\begin{center}$\begin{aligned}
\sup\{\Gpd(_RM)~|~_RM\in R\text{-Mod}\}=&\Glgl(R)=\sup\{\Gid(_RM)~|~_RM\in R\text{-Mod}\}\\
(\text{resp.,} \sup\{\Gfd(_RM)~|~_RM\in R\text{-Mod}\}=&\Gwgl(R)=\sup\{\Gfd(M_R)~|~M_R\in \text{Mod-}R\}).
\end{aligned}$\end{center}}
\end{df}

For a ring $R$, the invariant $\mathrm{sfli}(R)$  is defined as the supremum of the  flat dimensions of all injective left $R$-modules.

\begin{rem} \label{06-3-1}{\rm
(1) The equality
$\sup\{\Gpd(_RM)~|~_RM\in R\text{-Mod}\} = \sup\{\Gid(_R M)~|~_RM\in R\text{-Mod}\}$
was proved in \cite{Bel2000, BM2010, Emm2012} by using different methods.

(2) For any ring $R$, the equality
\begin{center}$\begin{aligned}
\sup\{\Gfd(_RM)~|~_RM \in R\text{-Mod}\}=\sup\{\Gfd(M_R~)~|~M_R\in \text{Mod-R}\}
\end{aligned}$\end{center}
was proved in~\cite{CET2021}.
Note from \cite{SS2020} that any ring is left and right GF closed (that is, the subcategories $_R \mathcal{GF}$ and $\mathcal{GF}_R$ are closed under extensions), so
the above equality can also be obtained from Bouchiba \cite[Theorem 6(2)]{Bou2015}, as noted in \cite{WZ2023}.
Furthermore, there is an equality
$$(\sharp)\quad \Gwgl(R)=\max\{\mathrm{sfli}(R), \mathrm{sfli}(R^{op})\}.$$

(3)
Note that there is an equality $$\sup\{\fd(_RM)~|~_RM\in~_R\mathcal{I}\}=\sup\{\fd(_RM)~|~_RM\in~_R\mathcal{FI}\}.$$
for any ring $R$.
 On one hand, the ``$\leq$'' is clear.
Conversely, let $_RM$  be any FP-injective left $R$-module.
Then there is a pure short exact sequence $0\to~_RM\to~ _RI\to ~_RK\to0$ of left $R$-modules with $_RI$ injective.
It follows from \cite[Lemma 9.1.4]{EJ2000} that $\fd(_RM)\leq\fd(_RI)$. Thus ``$\geq$'' holds true as well.

(4) We note that $\mathrm{sfli}(R)<\infty$ if and only if every injective left $R$-module has finite flat dimension. Indeed, the ``only if'' part is trivial. For the ``if'' part,
we suppose that every injective left $R$-module has finite flat dimension and assume towards a contradiction that $\sup\{\fd(_RM)~|~_RM\in~_R\mathcal{I}\}=\infty$ holds.
Then, by (3), for any positive integer $k$, there exists an FP-injective left $R$-module $_RM_k$
such that $\fd(_RM_k)\geq k$.
So one gets $\fd(\coprod_{k\in \mathbb{N}} M_k)=\sup\{\fd({_R} M_k)~|~k\in \mathbb{N}\}=\infty$.
Note that the subcategory $_R\mathcal{FI}$ is closed under arbitrary coproducts. In particular, $\coprod_{k\in \mathbb{N}}M_k \in ~_R\mathcal{FI}$, and so, as it is a pure submodule of an injective module, the assumption yields that it has finite flat dimension (again by \cite[Lemma 9.1.4]{EJ2000}), which is a contradiction.
Thus, $\sup\{\fd(_RM)~|~_RM\in~_R\mathcal{I}\}<\infty$ and hence $\mathrm{sfli}(R)<\infty$.}
\end{rem}


\subsection{\bf dg-Gorenstein homological complexes}

Gillespie \cite{Gil2004} introduced the following definitions,
which extend the notions of
dg-projective and dg-injective complexes and dg-projective and dg-injective resolutions of complexes, respectively.

\begin{df}\label{10-9} {\rm \cite[Definition 3.3]{Gil2004}. Let $(_R\mathcal{X}, ~_R\mathcal{Y})$ be a cotorsion pair in $R$-$\Mod$ and $_RM^{\bullet}$ a complex in \Ch$(R\text{-}\Mod)$.
\begin{enumerate}
\item $_RM^{\bullet}$ is called an \emph{$_R\mathcal{X}$ complex} if it is exact and $\Z^{k}(_RM^{\bullet}) \in ~_R\mathcal{X}$ for each $k \in \mathbb{Z}.$
\item $_RM^{\bullet}$ is called an \emph{$_R\mathcal{Y}$ complex} if it is exact and $\Z^{k}(_RM^{\bullet}) \in ~_R\mathcal{Y}$ for each $k \in \mathbb{Z}.$
\item $_RM^{\bullet}$ is called a \emph{dg-$\mathcal{X}$ complex} if $_RM^{k} \in ~_R\mathcal{X}$ for each $k \in \mathbb{Z}$, and
every morphism $f^{\bullet}:~_RM^{\bullet}\to ~_RY^{\bullet}$ with $_RY^{\bullet}$ an $_R\mathcal{Y}$ complex is null-homotopic.
\item $_RM^{\bullet}$ is called a \emph{dg-$\mathcal{Y}$ complex} if $_RM^{k} \in ~_R\mathcal{Y}$ for each $k \in \mathbb{Z}$, and
every morphism $f^{\bullet}:~_RX^{\bullet}\to ~_RM^{\bullet}$ with $_RX^{\bullet}$ an $_R\mathcal{X}$ complex is null-homotopic.
\item The subcategory of all $_R\mathcal{Y}$ (resp., \dg-$_R\mathcal{Y}$) complexes is denoted by $\widetilde{_R\mathcal{Y}}$ (resp., \dg$(_R\mathcal{Y})$).
Similarly, the subcategory of all $_R\mathcal{X}$ (resp., \dg-$_R\mathcal{X}$) complexes is denoted by $\widetilde{_R\mathcal{X}}$ (resp., \dg$(_R\mathcal{X})$).
\item A $dg$-$_R\mathcal{X}$ (resp., $dg$-$_R\mathcal{Y}$) resolution of $_RM^{\bullet}$ is a quasi-isomorphism $_RX^{\bullet} \to ~_RM^{\bullet}$ (resp. $_RM^{\bullet} \to ~_RY^{\bullet}$) in which $_RX^{\bullet}\in \dg(_R\mathcal{X})$ (resp., $_RY^{\bullet}\in \dg(_R\mathcal{Y})$).
\end{enumerate}
}
\end{df}

If we take $(\mathcal{X},\mathcal{Y})$ as $(_R\mathcal{P}, R\text{-Mod})$ (resp. $(R\text{-Mod},~_R\mathcal{I})$ and $(_R\mathcal{F}, _R\mathcal{C}ot)$),
then (1) dg-$_R\mathcal{P}$ (resp., dg-$_R\mathcal{I}$) complexes are exactly \emph{dg-projective} (resp., \emph{dg-injective}) complexes; (2) $\widetilde{_R\mathcal{P}}$ (resp., $\widetilde{_R\mathcal{I}}$ and $\widetilde{_R\mathcal{F}}$) complexes are exactly projective (resp., injective and flat) complexes, note that projective (resp. injective) complexes are \emph{contractible}, that is,
an exact complex such that the identity morphism is null-homotopic;
and (3) a dg-$_R\mathcal{P}$ (resp., dg-$_R\mathcal{I}$ and dg-$_R\mathcal{F}$) resolution of a complex $_RM^{\bullet}$ is exactly a dg-projective (resp., dg-injective and dg-flat) resolution of $_RM^{\bullet}$.

Next we consider the Gorenstein analogue notions.

\begin{df}\label{10-9-1}
{\rm
(1) According to \cite[Theorem 5.6]{SS2020}, for any ring $R$, there is a complete and hereditary cotorsion pair $(^{\perp}{_R\mathcal{GI}},~_R\mathcal{GI})$ in $R$-Mod.
So one has the notion of dg-$_R\mathcal{GI}$ complexes.
We say that a complex $_RG^{\bullet}$  is \emph{dg}-\emph{Gorenstein injective} if $G^{\bullet}$ is  dg-$_R\mathcal{GI}$.
In addition, a morphism $_RX^{\bullet} \to ~_RG^{\bullet}$ is called \emph{dg}-\emph{Gorenstein injective resolution}
of the complex $_RX^{\bullet}$ if it is a dg-$_R\mathcal{GI}$ resolution of $_RX^{\bullet}$,
that is, $_RX^{\bullet}\to ~_RG^{\bullet}$ is a quasi-isomorphism with $_RG^{\bullet}$ dg-Gorenstein injective.

(2) Let $R$ be a ring such that there is a complete and hereditary cotorsion pair $(_R\mathcal{GP},~{_R\mathcal{GP}^{\perp}})$ in $R$-Mod. Then one gets the notion of dg-$_R\mathcal{GP}$ complexes.
We say that a complex $_RG^{\bullet}$ is \emph{dg}-\emph{Gorenstein projective} if $_RG^{\bullet}$ is  dg-$_R\mathcal{GP}$.
In addition, a morphism $_RG^{\bullet} \to ~_RX^{\bullet}$ is called \emph{dg}-\emph{Gorenstein projective resolution}
of the complex $_RX^{\bullet}$  if it is a dg-$_R\mathcal{GP}$ resolution of $_RX^{\bullet}$,
that is, $_RG^{\bullet}\to ~_RX^{\bullet}$ is a quasi-isomorphism with $_RG^{\bullet}$ dg-Gorenstein projective.
}
\end{df}

\subsection{\bf Gorenstein derived categories}
As mentioned above, we denote by $\mathbf{D}(R\text{-}\Mod)$ the derived category of all left $R$-modules.
Gao and Zhang \cite{GZ2010} introduced the notion of Gorenstein derived categories.
These notions can be unified as the following:

\begin{df}\label{10-9-2}
{\rm Let $R$ be a ring and $_R\mathcal{X}$, $_R\mathcal{Y}$ be subcategories of $R$-Mod which contain $_R\mathcal{P}$ and $_R\mathcal{I}$ respectively.
\begin{enumerate}
\item A complex $_RC^{\bullet}$ in \Ch$(R\text{-Mod})$ is called \emph{right $_R\mathcal{X}$-acyclic} (resp. \emph{left $_R\mathcal{Y}$-acyclic}) if $_RC^{\bullet}$ is $\Hom_R(_R\mathcal{X},-)$-exact. (resp.~$\Hom_R(-,~_R\mathcal{Y})$-exact).
\item A morphism~$f^{\bullet} : ~_RC^{\bullet}\rightarrow ~_RD^{\bullet}$ in \Ch$(R\text{-Mod})$ is said to be a \emph{right $_R\mathcal{X}$-quasi-isomorphism}~(resp. \emph{left $_R\mathcal{Y}$-quasi-isomorphism})
if Hom$_{R}(_RX, f^{\bullet})$ is a quasi-isomorphism for all $_RX\in~_R\mathcal{X}$ (resp.  Hom$_R(f^{\bullet}, ~_RY)$ is a quasi-isomorphism for all $_RY\in~_R\mathcal{Y}$).
\item Denote by $_R\mathcal{X}\text{-}\mathrm{rac}$ (resp. $_R\mathcal{Y}\text{-}\mathrm{lac}$) the subcategory of $\Ch(R\text{-}\Mod)$
consisting of all right $_R\mathcal{X}$-acyclic (resp. left $_R\mathcal{Y}$-acyclic) complexes. So $\K(_R\mathcal{X}\text{-}\mathrm{rac})$ and $\K(_R\mathcal{X}\text{-}\mathrm{rac})$
mean the the corresponding homotopy subcategories respectively. Note that such homotopy subcategories are both a thick subcategory of $\K(R\text{-}\Mod)$.
Now, the relative derived category with respect to $_R\mathcal{X}$ (resp. $_R\mathcal{Y}$),
denoted by $\mathbf{D}_{\mathcal{X}}(R\text{-}\Mod)$ (resp. $\mathbf{D}_{\mathcal{Y}}(R\text{-}\Mod)$), is defined as the following Verdier quotient triangulated category:
$$\mathbf{D}_{\mathcal{X}}(R\text{-}\Mod)=\K(R\text{-}\Mod)/\K(\mathcal{X}\text{-}\mathrm{rac})~(\mathrm{resp}. \quad\mathbf{D}_{\mathcal{Y}}(R\text{-}\Mod)=\K(R\text{-}\Mod)/\K(\mathcal{Y}\text{-}\mathrm{lac})).$$
\end{enumerate}
 }
\end{df}

If we take the pair $(_R\mathcal{X},~_R\mathcal{Y})$ as $(_R\mathcal{P},~_R\mathcal{I})$, then the relative derived categories with respect to $_R\mathcal{X}$ and $_R\mathcal{Y}$,
coincide with each other, which are exactly the usual derived category of $R$-Mod, $\mathbf{D}(R\text{-}\Mod)$.

If we take the pair $(\mathcal{X},\mathcal{Y})$ as $(_R\mathcal{GP},~_R\mathcal{GI})$, then the relative derived categories with respect to $_R\mathcal{X}$ and $_R\mathcal{Y}$ respectively, are called \emph{Gorenstein projective} and \emph{Gorenstein injective} \emph{derived category} of $R$-Mod, respectively.
We denote such relative derived categories by $\mathbf{D}_{\mathrm{GP}}(R\text{-}\Mod)$ and $\mathbf{D}_{\mathrm{GI}}(R\text{-}\Mod)$ respectively.
As noted in \cite{GZ2010}, in contrast to the classic situation, $\mathbf{D}_{\mathrm{GP}}(R\text{-}\Mod)$ and $\mathbf{D}_{\mathrm{GI}}(R\text{-}\Mod)$ are not the same, in general.

\section{\bf Realizing relative derived categories as homotopy categories.}
As mentioned above, we denote by $\dg(_R\mathcal{GP})$ (resp. $\dg(_R\mathcal{I})$) the subcategory of \Ch$(R\text{-}\Mod)$ consisting of all
dg-projective (resp. dg-injective) complexes. Hence, $\K(\dg(_R\mathcal{P}))$ (resp. $\K(\dg(_R\mathcal{P}))$) means the homotopy subcategory
consisting of all dg-projective (resp. dg-injective) complexes.

It is well known that there are triangle equivalences
$$\K(\dg(_R\mathcal{P}))\simeq\mathbf{D}(R\text{-}\Mod)\simeq\K(\dg(_R\mathcal{I}))$$
for any ring $R$.
This realizes the usual derived categories as homotopy categories.
In this section, we will consider the result for relative derived categories.
Our aim is to obtain some torsion pairs.

\begin{df} \label{T-0}\rm{
Let $\mathcal{T}$ be a triangulated category and $\mathcal{T}_1, \mathcal{T}_2$ its triangulated subcategory. The pair $(\mathcal{T}_1, \mathcal{T}_2)$ is said to be a \emph{torsion pair}
if the following are satisfied
\begin{enumerate}
  \item[(TP1)] $\Hom_{\mathcal{T}}(T_1,T_2)=0$ for all $T_1\in \mathcal{T}_1$ and all $T_2\in \mathcal{T}_2$.
  \item[(TP2)] For any $T\in \mathcal{T}$, there is a distinguished triangle
  $T_1\to T\to T_2\to\Sigma T_1$
  with $T_1\in \mathcal{T}_1$ and $T_2\in \mathcal{T}_2$, where $\Sigma$ denotes the autofunctor of $\mathcal{T}$.
\end{enumerate}}
\end{df}

\begin{rem} \label{T-0-0}\rm{
According to \cite[Chapter I, Definitions 2.1 and 2.5, and Proposition 2.6]{BI2007}, the notion in Definition \ref{T-0} should be named as a ``hereditary torsion pair''.
We simply call it a ``torsion pair''.}
\end{rem}

As shown in the following lemma, torsion pairs are an important tool to obtain some triangle-equivalences.

\begin{lem} \label{T-0-1}
Let $(\mathcal{T}_1, \mathcal{T}_2)$ be a torsion pair of a triangulated category $\mathcal{T}$. Then the following
composite functors
$$F_1:\mathcal{T}_1\overset{\mathrm{inc}_1}\to \mathcal{T}\overset{Q_2}\to \mathcal{T}/\mathcal{T}_2\quad\text{and}\quad F_2:\mathcal{T}_2\overset{\mathrm{inc}_2}\to \mathcal{T}\overset{Q_1}\to \mathcal{T}/\mathcal{T}_1$$
are triangle equivalent, where $\mathrm{inc}_1$ and $\mathrm{inc}_2$ mean the canonical inclusion functors respectively, and $Q_1$ and $Q_2$ denote the canonical quotient functors respectively.
\end{lem}

\begin{proof} The result is well-known. For the left-hand triangle equivalence,
one can see the proof of (i)$\Rightarrow$(vi) in \cite[Chapter I, Proposition 2.6]{BI2007}, and the right-hand one is a dual.
\end{proof}

\begin{rem} \label{T-0-2}{\rm
Let $(\mathcal{T}_1, \mathcal{T}_2)$ be a torsion pair of a triangulated category $\mathcal{T}$
and $F_1$ and $F_2$ be the equivalent functors in Lemma \ref{T-0-1}. We write $F'_1$ and $F'_2$ as their quasi-inverses respectively.
According to the proof of (i)$\Rightarrow$(vi) and its dual in \cite[Chapter I, Proposition 2.6]{BI2007}, we know that the rule of $F'_1$ and $F'_2$ for objects are as follows:
$$F'_1(Y)=T_2^{Y}, \forall Y\in  \mathcal{T}/\mathcal{T}_2\quad\text{and}\quad F'_2(X)=T_2^{X}, \forall X\in  \mathcal{T}/\mathcal{T}_1,$$
where $T_2^{Y}$ and $T_1^{X}$ come from the distinguished triangles
$$T^{Y}_1\to Y\to T_2^{Y}\to\Sigma T^{Y}_1\quad\text{and}\quad T^{X}_1\to X\to T^{X}_2\to\Sigma T^{X}_1$$
which are guaranteed by (TP2).
}
\end{rem}

In order to obtain our desired torsion pairs, we need the following lemmas.

\begin{lem} \label{T-1} Let $R$ be a ring and $_R\mathcal{X}$ a subcategory of $R\text{-}\Mod$.
If there is a complete and hereditary cotorsion pair $(_R\mathcal{X}, ~_R\mathcal{X}^{\perp})$ in
$R\text{-}\Mod$, then every complex in $\widetilde{_R\mathcal{X}^{\perp}}$ is right $_R\mathcal{X}$-acyclic.
\end{lem}
\begin{proof} Let $$_RX^{\bullet}=\cdots \to ~_RX^{k-1}\overset{d^{k-1}}\to ~_RX^{k}\overset{d^{k}}\to ~_RX^{k+1}\to\cdots$$ be a complex in $\widetilde{_R\mathcal{X}^{\perp}}$.
Then by the definition, $_RX^{\bullet}$ is exact and each cycle module $\Ker{d^{k}}$ is in $_R\mathcal{X}^{\perp}$. Now we decompose $_RX^{\bullet}$ into short exact sequences $$_RX^{\bullet}_k = 0\to \Ker{d^{k}}\to ~_RX^{k}\overset{d^{k}}\to \Ker{d^{k+1}}\to 0, k\in \mathbb{Z}.$$
Note that $\Ext_R^{1}(_RX',\Ker{d^{k}})=0$ for all $_RX'\in ~_R\mathcal{X}$ and all $k\in \mathbb{Z}$ as each $\Ker{d^{K}} \in ~_R\mathcal{X}^{\perp}$.
This implies that each $_RX^{\bullet}_k$ is Hom$_R(_R\mathcal{X},-)$-exact,
and then so is $_RX^{\bullet}$. It follows that $_RX^{\bullet}$ is right $_R\mathcal{X}$-acyclic.
\end{proof}

\begin{lem} \label{T-2}  Let $R$ be a ring and $_R\mathcal{X}$ a subcategory of $R\text{-}\Mod$.
If there is a complete and hereditary cotorsion pair $(_R\mathcal{X}, ~_R\mathcal{X}^{\perp})$ in
$R\text{-}\Mod$, then every complex in $\dg(_R\mathcal{X})\cap ~_R\mathcal{X}\text{-}\mathrm{rac}$ is contractible.
\end{lem}

\begin{proof}
We let $$_RX^{\bullet}=\cdots \to ~_RX^{k-1}\overset{d^{k-1}}\to ~_RX^{k}\overset{d^{k}}\to ~_RX^{k+1}\to\cdots$$
be a complex in $\dg(_R\mathcal{X})\cap ~_R\mathcal{X}\text{-}\mathrm{rac}$.
For each $k\in \mathbb{Z}$, let us consider the complex $$0\to \Ker{d^{k}}\to ~_RX^{k}\overset{d^{k}}\to \Ker{d^{k+1}}\to 0.$$
Note that any complex in $_R\mathcal{X}\text{-rac}$ is exact as $_R\mathcal{P}\subseteq ~_R\mathcal{X}$,
and hence, by \cite[Theorem 3.12]{Gil2004} one has $$_RX^{\bullet}\in \dg(_R\mathcal{X})\cap ~_R\mathcal{X}\text{-rac} \subseteq \widetilde{_R\mathcal{X}}\cap ~_R\mathcal{X}\text{-rac}.$$
It follows that, for each $k\in \mathbb{Z}$, the complex $0\to \Ker{d^{k}}\to ~_RX^{k}\overset{d^{k}}\to \Ker{d^{k+1}}\to 0$ is a short exact sequence in $R\text{-}\Mod$ with $\Ker{d^{k+1}}\in~_R\mathcal{X}$.
Thus, the identity morphism of $\Ker{d^{k+1}}$ factors through $d^{k}$. That is, $_RX^{\bullet}$ is an exact complex with each $d^{k}$ a split epic morphism, and so, $_RX^{\bullet}$ is contractible.
\end{proof}

\begin{lem} \label{T-3}  Let $R$ be a ring and $_R\mathcal{X}$ a subcategory of $R\text{-}\Mod$.
If there is a complete and hereditary cotorsion pair $(_R\mathcal{X}, ~_R\mathcal{X}^{\perp})$ in
$R\text{-}\Mod$, then $\Hom_{\mathrm{K}(R\text{-}\mathrm{Mod})}(_RX^{\bullet},~_RY^{\bullet})=0$ for all $_RX^{\bullet}\in \dg(_R\mathcal{X})$ and all $_RY^{\bullet}\in ~_R\mathcal{X}\text{-}\mathrm{rac}$.
\end{lem}

\begin{proof}  Let $_RX^{\bullet}\in \dg(_R\mathcal{X})$, $_RY^{\bullet}\in ~_R\mathcal{X}\text{-}\mathrm{rac}$ and consider a morphism $f^{\bullet}: ~_RX^{\bullet}\to ~_RY^{\bullet}$.
Since $(_R\mathcal{X}, ~_R\mathcal{X}^{\perp})$ is a complete and hereditary cotorsion pair in $R\text{-}\Mod$, it,
applying \cite[Theorem 3.5]{YgL2011}, induces a complete and hereditary cotorsion pair $(\dg(_R\mathcal{X}), \widetilde{_R\mathcal{X}^{\perp}})$ in  $\Ch(R\text{-}\Mod)$.
In particular, there is a short exact sequence in $\Ch(R\text{-}\Mod)$
$$0\to~ _RK^{\bullet} \to ~_RG^{\bullet}\overset{\alpha^{\bullet}}\to ~_RY^{\bullet}\to 0$$
with $_RG^{\bullet}\in \dg(_R\mathcal{X})$ and $_RK^{\bullet}\in \widetilde{_R\mathcal{X}^{\perp}}$.
Note that $\alpha^{\bullet}: ~_RG^{\bullet}\to _RY^{\bullet}$ is a special $\dg(_R\mathcal{X})$-precover of $_RY^{\bullet}$,
there is a morphism $\beta^{\bullet}:~ _RX^{\bullet}\to ~_RG^{\bullet}$ such that $f^{\bullet}=\alpha^{\bullet}\beta^{\bullet}$.

By Lemma \ref{T-1} one gets $_RK^{\bullet}\in ~_R\mathcal{X}\text{-}\mathrm{rac}$. It follows that $_RG^{\bullet}$ is in $_R\mathcal{X}\text{-}\mathrm{rac}$ since so is $_RY^{\bullet}$, and hence, $_RG^{\bullet}\in \dg(_R\mathcal{X})\cap ~_R\mathcal{X}\text{-rac}$.
Thus, Lemma \ref{T-2} yields that $_RG^{\bullet}$ is is contractible.
Consequently, the morphism $\alpha^{\bullet}$ is null-homotopic and so is $f^{\bullet}=\alpha^{\bullet}\beta^{\bullet}$. That is, $\Hom_{\mathrm{K}(R\text{-}\mathrm{Mod})}(_RX^{\bullet},~_RY^{\bullet})=0$.
\end{proof}

The next three results are the dual of Lemmas \ref{T-1}-\ref{T-3} respectively.

\begin{lem} \label{T-1'} Let $R$ be a ring and $_R\mathcal{Y}$ a subcategory of $R\text{-}\Mod$.
If there is a complete and hereditary cotorsion pair  $(^{\perp}{_R\mathcal{Y}}, ~_R\mathcal{Y})$ in
$R\text{-}\Mod$, then every complex in $\widetilde{^{\perp}{_R\mathcal{Y}}}$ is left $_R\mathcal{Y}$-acyclic.
\end{lem}

\begin{lem} \label{T-2'}  Let $R$ be a ring and $_R\mathcal{Y}$ a subcategory of $R\text{-}\Mod$.
If there is a complete and hereditary cotorsion pair  $(^{\perp}{_R\mathcal{Y}}, ~_R\mathcal{Y})$ in
$R\text{-}\Mod$, then every complex in $_R\mathcal{Y}\text{-}\mathrm{lac}\cap \dg(_R\mathcal{Y})$ is contractible.
\end{lem}

\begin{lem} \label{T-3'}  Let $R$ be a ring and $_R\mathcal{Y}$ a subcategory of $R\text{-}\Mod$.
If there is a complete and hereditary cotorsion pair  $(^{\perp}{_R\mathcal{Y}}, ~_R\mathcal{Y})$ in
$R\text{-}\Mod$, then $\Hom_{\mathrm{K}(R\text{-}\mathrm{Mod})}(_RX^{\bullet},~_RY^{\bullet})=0$ for all $_RX^{\bullet}\in~ _R\mathcal{Y}\text{-}\mathrm{lac}$ and all $_RY^{\bullet}\in \dg(_R\mathcal{Y})$.
\end{lem}

Now we can give some torsion pairs of $\K(R\text{-}\Mod)$.

\begin{prop} \label{T-4} Let $R$ be a ring and $_R\mathcal{X}$ a subcategory of $R\text{-}\Mod$.
If there is a complete and hereditary cotorsion pair $(_R\mathcal{X}, ~_R\mathcal{X}^{\perp})$  in
$R\text{-}\Mod$, then there is a torsion pair $(\K(\dg(_R\mathcal{X})), \K(_R\mathcal{X}\text{-}\mathrm{rac}))$
in $\K(R\text{-}\Mod)$.
\end{prop}
\begin{proof} We must check that the pair $(\K(\dg(_R\mathcal{X})),\K(_R\mathcal{X}\text{-rac}))$ satisfies the conditions (TP1) and (TP2) in Definition \ref{T-0}.
(TP1) holds by Lemma \ref{T-3}. To see (TP2), let $_RM^{\bullet}\in $ $\K(R\text{-}\Mod)$.
As mentioned above, \cite[Theorem 3.5]{YgL2011} allows us to have  a short exact sequence in $\Ch(R\text{-}\Mod)$
$$0\to ~_RK^{\bullet} \to ~_RG^{\bullet}\overset{\alpha^{\bullet}}\to ~_RM^{\bullet}\to 0$$
with $_RG^{\bullet}\in \dg(_R\mathcal{X})$ and $_RK^{\bullet}\in \widetilde{_R\mathcal{X}^{\perp}}$.
Now we consider the distinguished triangle in $\K(R\text{-}\Mod)$
$$_RG^{\bullet}\overset{\alpha^{\bullet}}\to ~_RM^{\bullet}\to \mathrm{Cone}(\alpha^{\bullet})\to \Sigma _RG^{\bullet}$$
where $_RG^{\bullet}$ is in $\dg(_R\mathcal{X})$.
For any left $R$-module $_RX\in ~_R\mathcal{X}$, viewed it as a complex, it is in$~\dg(_R\mathcal{X})$, and hence, $\Hom_R(_RX, ~_RK^{\bullet})=\mathcal{H}om(_RX, ~_RK^{\bullet})$
is exact. It follows that $\Hom_R(_RX,\alpha^{\bullet})$ is a quasi-isomorphism for all $_RX\in ~_R\mathcal{X}$, and so $\alpha^{\bullet}$ is a right $_R\mathcal{X}$-quasi-isomorphism.
This implies that $\mathrm{Cone}(\alpha^{\bullet}) \in \K(_R\mathcal{X}\text{-}\mathrm{rac})$ and then (TP2) follows.
\end{proof}

\begin{prop} \label{T-4'} Let $R$ be a ring and $_R\mathcal{Y}$ a subcategory of $R\text{-}\Mod$.
If there is a complete and hereditary cotorsion pair  $(^{\perp}{_R\mathcal{Y}}, ~_R\mathcal{Y})$ in
$R\text{-}\Mod$, then there is a torsion pair $(\K(_R\mathcal{Y}\text{-}\mathrm{lac}),\K(\dg(_R\mathcal{Y})))$
in $\K(R\text{-}\Mod)$.
\end{prop}
\begin{proof} By an argument dual to Proposition \ref{T-4}.
\end{proof}

Note that $\K(R\text{-}\Mod)/\K(\mathcal{X}\text{-}\mathrm{rac})= \mathbf{D}_{\mathcal{X}}(R\text{-}\Mod)$ and $\K(R\text{-}\Mod)/\K(\mathcal{Y}\text{-}\mathrm{lac})= \mathbf{D}_{\mathcal{Y}}(R\text{-}\Mod)$.
Now applying Lemma \ref{T-0-1} to the torsion pairs from Propositions \ref{T-4} and \ref{T-4'}, one can obtain the following result, which realizes the relative derived category $\mathbf{D}_{\mathcal{X}}(R\text{-}\Mod)$ (resp. $\mathbf{D}_{\mathcal{Y}}(R\text{-}\Mod)$) as a homotopy subcategory.

\begin{thm} \label{T-5} Let $R$ be a ring and $_R\mathcal{X}, ~_R\mathcal{Y}$ subcategories of $R\text{-}\Mod$.
\begin{enumerate}
  \item If there is a complete and hereditary cotorsion pair $(_R\mathcal{X}, ~_R\mathcal{X}^{\perp})$ in $R\text{-}\Mod$, then there is a triangle-equivalence $\K(\dg(_R\mathcal{X}))\simeq\K(R\text{-}\Mod)/\K(_R\mathcal{X}\text{-}\mathrm{rac})=\mathbf{D}_{\mathcal{X}}(R\text{-}\Mod)$.
  \item If there is a complete and hereditary cotorsion pair $(^{\perp}{_R\mathcal{Y}}, ~_R\mathcal{Y})$ in $R\text{-}\Mod$, then there is a triangle-equivalence $\K(\dg(_R\mathcal{Y}))\simeq\K(R\text{-}\Mod)/\K(_R\mathcal{Y}\text{-}\mathrm{lac})=\mathbf{D}_{\mathcal{Y}}(R\text{-}\Mod)$.
\end{enumerate}
\end{thm}

\begin{cor} \label{T}
 Let $R$ be a ring and $_R\mathcal{X}, ~_R\mathcal{Y}$ subcategories of $R\text{-}\Mod$.
If there is a complete and hereditary cotorsion triple $(_R\mathcal{X}, ~_R\mathcal{X}^{\perp}=~^{\perp}{_R\mathcal{Y}}, ~_R\mathcal{Y})$ in
$R\text{-}\Mod$, then there are triangle-equivalences
$$\K(\dg(_R\mathcal{X}))\simeq\mathbf{D}_{\mathcal{X}}(R\text{-}\Mod)=\mathbf{D}_{\mathcal{Y}}(R\text{-}\Mod)\simeq \K(\dg(_R\mathcal{Y})).$$
\end{cor}
\begin{proof}
Note that $\K(R\text{-}\Mod)/\K(_R\mathcal{X}\text{-}\mathrm{rac})= \mathbf{D}_{\mathcal{X}}(R\text{-}\Mod)$ and $\K(R\text{-}\Mod)/\K(_R\mathcal{Y}\text{-}\mathrm{lac})= \mathbf{D}_{\mathcal{Y}}(R\text{-}\Mod)$.
Since $(_R\mathcal{X}, ~_R\mathcal{X}^{\perp}=~^{\perp}{_R\mathcal{Y}}, ~_R\mathcal{Y})$ is a complete and hereditary cotorsion triple in $R\text{-}\Mod$, one has $\K(_R\mathcal{X}\text{-}\mathrm{rac})=\K(_R\mathcal{Y}\text{-}\mathrm{lac})$
by  \cite[Propositions 2.6 and 2.2]{Chen2010}, and hence $\mathbf{D}_{\mathcal{X}}(R\text{-}\Mod)=\mathbf{D}_{\mathcal{Y}}(R\text{-}\Mod)$.
At the same time, the triangle-equivalences $\K(\dg(_R\mathcal{X}))\simeq\mathbf{D}_{\mathcal{X}}(R\text{-}\Mod)$ and $\K(\dg(_R\mathcal{Y}))\simeq\mathbf{D}_{\mathcal{Y}}(R\text{-}\Mod)$ holds by Theorem \ref{T-5}.
\end{proof}

\begin{rem}\label{T-7}
{\rm
Let $R$ be a ring such that $(_R\mathcal{X}, ~_R\mathcal{X}^{\perp}=~^{\perp}{_R\mathcal{Y}}, ~_R\mathcal{Y})$ is a complete and hereditary cotorsion triple in
$R\text{-}\Mod$. Then Corollary \ref{T} shows a triangle-equivalence $\K(\dg(_R\mathcal{X}))\simeq \K(\dg(_R\mathcal{Y}))$.
We denote by $H$ and $H'$ the triangle-equivalent functor $\K(\dg(_R\mathcal{X}))\to \K(\dg(_R\mathcal{Y}))$ and its quasi-inverse.

According to Lemma \ref{T-0-1} and Remark \ref{T-0-2}, $H$ and $H'$ can be be described as follows:
\begin{center}$\begin{aligned}& H: \K(\dg(_R\mathcal{X}))\overset{\mathrm{inc}_{\mathcal{X}}}\to \K(R\text{-}\Mod)\overset{Q_{\mathcal{X}}}\to \K(R\text{-}\Mod)/\K(_R\mathcal{X}\text{-}\mathrm{rac})= \K(R\text{-}\Mod)/\K(_R\mathcal{Y}\text{-}\mathrm{lac})\overset{F'_{\mathcal{Y}}}\rightarrow \K(\dg(_R\mathcal{Y}))\\
& H': \K(\dg(_R\mathcal{Y}))\overset{\mathrm{inc}_{\mathcal{Y}}}\to \K(R\text{-}\Mod)\overset{Q_{\mathcal{Y}}}\to \K(R\text{-}\Mod)/\K(_R\mathcal{Y}\text{-}\mathrm{lac})= \K(R\text{-}\Mod)/\K(_R\mathcal{X}\text{-}\mathrm{rac})\overset{F'_{\mathcal{X}}}\rightarrow \K(\dg(_R\mathcal{X})).
\end{aligned}$\end{center}

Hence, the rule of $H$ and $H'$ for objects are as follows
\begin{center}$\begin{aligned}&H: ~_RA^{\bullet}\mapsto ~_RB^{\bullet}, \text{where}~_RB^{\bullet}~\text{is a dg-}{_R\mathcal{Y}}~\text{resolution of}~_RA^{\bullet}; \\
&H': ~_RN^{\bullet}\mapsto ~_RM^{\bullet}, \text{where}~_RM^{\bullet}~~\text{is a dg-}{_R\mathcal{X}}~\text{resolution of}~_RN^{\bullet}.
\end{aligned}$\end{center}}
\end{rem}

Recall from \cite{DLW2023} that a ring $R$ is \emph{left virtually Gorenstein} if $~_R\mathcal{GP}^{\perp}=~^{\perp}{_R\mathcal{GI}}$.
We denote by $_R\mathcal{PGF}$ the subcategory of $R$-Mod consisting of all projectively coresolved Gorenstein flat modules.
Here, a left $R$-module $_RM$ is said to be a \emph{ projectively coresolved Gorenstein flat} \cite{SS2020}
if there exists an \textbf{F}-totally acyclic complex $_RP^{\bullet}$ of projective left $R$-modules
such that $_RM\cong Z^{0}(_RP^{\bullet})$.

Let $_R\mathcal{X}$ be a subcategory of $R\text{-}\Mod$.
Following \cite{SS2020},
the \emph{definable closure} of $_R\mathcal{X}$, denoted by $\overline{_R\mathcal{X}}$,
is the smallest subcategory which is definable and contains $_R\mathcal{X}$.
Here we say that $\overline{_R\mathcal{X}}$ is~\emph{definable} if $\overline{_R\mathcal{X}}$ is closed under pure submodules, direct limits and direct products.

\begin{prop}\label{T-8-0}
 Let $R$ be a left virtually Gorenstein ring. Then we have the equality $_R\mathcal{GP}=~_R\mathcal{PGF}$. In particular, there is a projective cotorsion pair $(_R\mathcal{GP}, ~{_R\mathcal{GP}^{\perp}}).$
\end{prop}
\begin{proof} Assume that $R$ is a left virtually Gorenstein ring.
Note from \cite[Theorem 5.6 and its proof]{SS2020} that, for any ring $R$,  the subcategory $^{\perp}{_R\mathcal{GI}}$ is closed under direct limits.
Therefore, as $R$ is left virtually Gorenstein, the subcategory ${_R\mathcal{GP}^{\perp}}$ is closed under direct limits.
On the other hand,  we know from \cite[Corollary 3.4(1)]{CS2021} that the pair $(_R\mathcal{GP}, ~{_R\mathcal{GP}^{\perp}})$ forms a cotorsion pair (over any ring).
Thus, \cite[Theorem 6.1]{Sar2018} gives that the subcategory ${_R\mathcal{GP}^{\perp}}$ is definable. In particular, one gets $\overline{_RR} \subseteq ~{_R\mathcal{GP}^{\perp}}$.
So $_R\mathcal{GP}=~_R\mathcal{PGF}$ by \cite[Lemma 3.3]{WYZ2023}, and hence the cotorsion pair
$$(_R\mathcal{GP}, ~{_R\mathcal{GP}^{\perp}})=(_R\mathcal{PGF}, ~{_R\mathcal{PGF}^{\perp}})$$
is complete and hereditary by \cite[Theorem 4.9]{SS2020}. It is further projective by \cite[Theorem 5.4]{Gil2016}.
\end{proof}

\begin{rem}\label{T-8-1}
{\rm
A classical fact in algebra is that projective modules are always flat. However, it is still unknown
whether or not ``all Gorenstein projective left $R$-modules are Gorenstein flat'' holds over an arbitrary ring $R$.
According to Proposition \ref{T-8-0}, we know that all Gorenstein projective left $R$-modules are Gorenstein flat whenever $R$ is left virtually Gorenstein, since
projectively coresolved Gorenstein flat modules are always Gorenstein flat.
Thus, to ensure that the result holds, the statements ``Suppose that $R$ is a left virtually Gorenstein ring satisfying $_R\mathcal{GP}\subseteq ~_R\mathcal{GF}$'' and ``Suppose that $R$ is a virtually Gorenstein ring satisfying $_R\mathcal{GP}\subseteq ~_R\mathcal{GF}$ and $\mathcal{GP}_R\subseteq ~\mathcal{GF}_R$'' from \cite[Theorems B and C]{DLW2023} can be restated as ``Suppose that $R$ is a left virtually Gorenstein ring '' and ``Suppose that $R$ is a virtually Gorenstein ring '', respectively.
}
\end{rem}

We end this section with some examples.

\begin{exa}\label{T-9}
{\rm
(1) Let $R$ be any ring. There is a complete and hereditary triple cotorsion given by $(_R\mathcal{P}, R\text{-}\Mod, ~_R\mathcal{I})$ in
$R\text{-}\Mod$.
So Corollary \ref{T} recovers the well-known fact: the triangle-equivalences
$$\K(\dg(_R\mathcal{P}))\simeq\mathbf{D}(R\text{-}\Mod)\simeq\K(\dg(_R\mathcal{I}))$$
holds true for any ring $R$.

(2)
Let $R$ be a left virtually Gorenstein ring. Then there is a complete and hereditary cotorsion pair $(_R\mathcal{GP}, ~{_R\mathcal{GP}^{\perp}})$ in $R\text{-}\Mod$ by Proposition \ref{T-8-0}.
And hence, there is complete and hereditary cotorsion triple $(_R\mathcal{GP}, ~{_R\mathcal{GP}^{\perp}}={^{\perp}{_R\mathcal{GI}}},~ _R\mathcal{GI})$  in
$R\text{-}\Mod$, since, for any ring $R$, $(^{\perp}{_R\mathcal{GI}}, _R\mathcal{GI})$ forms a complete and hereditary cotorsion pair in
$R\text{-}\Mod$ (see \cite[Theorem 5.6]{SS2020}).
It follows from Corollary \ref{T} that there are triangle-equivalences
$$\K(\dg(_R\mathcal{GP}))\simeq\mathbf{D}_{\mathrm{GP}}(R\text{-}\Mod)\simeq\mathbf{D}_{\mathrm{GI}}(R\text{-}\Mod)\simeq\K(\dg(_R\mathcal{GI})).$$
Here $\dg(_R\mathcal{GP})$ (resp. $\dg(_R\mathcal{GI})$) denotes the subcategory of \Ch$(R\text{-}\Mod)$ consisting of all
dg-Gorenstein projective (resp. dg-Gorenstein injective) complexes. Hence, $\K(\dg(_R\mathcal{GP}))$ (resp. $\K(\dg(_R\mathcal{GI}))$) means the corresponding homotopy subcategory.
Furthermore, if we denote by $F:\K(\dg(_R\mathcal{GP}))\to \K(\dg(_R\mathcal{GI}))$ and $F':\K(\dg(_R\mathcal{GI}))\to \K(\dg(_R\mathcal{GP}))$ the associate equivalent functor respectively,
then the rules of $F$ and $F'$ on objects are defined as follows (due to Remark \ref{T-7}):
\begin{center}$\begin{aligned}&F: ~_RX^{\bullet}\mapsto ~_RY^{\bullet}, \text{where}~_RY^{\bullet}~\text{is a dg-Gorenstein injective resolution of}~_RX^{\bullet}; \\
&F': ~_RN^{\bullet}\mapsto _RM^{\bullet}, \text{where}~_RM^{\bullet}~\text{is a dg-Gorenstein projective resolution of}~_RN^{\bullet}.
\end{aligned}$\end{center}
}
\end{exa}

\section{\bf Applying to rings with finite Gorenstein weak global dimension}

According to \cite[Theorem A]{DLW2023}, we know that any ring with finite Gorenstein weak global dimension is both left and right virtually Gorenstein.
In this section, we will give some characterizations of a ring with finite Gorenstein weak global dimension and obtain some homotopic equivalences over such a ring.


In what follows, we denote by $\mathcal{F}^\wedge$ the subcategory consisting of all left or right $R$-modules with finite flat dimension
 and by $\mathcal{F}_m$ the subcategory consisting of all left or right $R$-modules with finite flat dimension at most $m$ for an integer $m\geq 0$;
by $\sfli(R)$  the supremum of the flat dimensions of all injective left $R$-modules;
by $\mathrm{FFD}(R)$ the left finitistic flat dimension of $R$, which is
defined as the supremum of the flat dimensions of those left $R$-modules that have finite flat dimension. The following result is immediate from \cite[Proposition 1.9]{C-W2023}, by taking $(\mathcal{U},\mathcal{V})=(_R\mathcal{F},_R\mathcal{F}^\perp)$, and using \cite[Theorem 4.1.3]{GT2006}. For the benefit of the reader, we opt for including a more independent argument below.

\begin{lem} \label{3-3-0}
Let $m$ be a nonnegative integer. Then the following are equivalent for any $R$:
\begin{enumerate}
  \item $\max\{\sfli(R), \mathrm{FFD}(R)\}\leq m$.
  \item The pair $(_R\mathcal{F}_m~, ~_R\mathcal{GI})$ is an injective cotorsion pair in $R\text{-}\Mod$.
  \item The pair $(_R\mathcal{F}_m~, ~{_R\mathcal{F}_m}^{\perp})$ is an injective cotorsion pair in $R\text{-}\Mod$.
\end{enumerate}
\end{lem}
\begin{proof}
(1)$\Rightarrow$(3). Assume $\max\{\sfli(R), \mathrm{FFD}(R)\}\leq m$. Due to \cite[Theorem 4.1.3]{GT2006}, the pair $(\mathcal{F}_m~,~{\mathcal{F}_m}^{\perp})$ is a complete hereditary cotorsion pair in $R$-Mod.
Thus, it is enough to prove $_R\mathcal{F}_m \cap ~{_R\mathcal{F}_m}^{\perp}=~_R\mathcal{I}$ by \cite[Proposition 3.6]{Gil2016}.
Note that $_R\mathcal{I} \subseteq ~_R\mathcal{F}_m \cap ~{_R\mathcal{F}_m}^{\perp}$
as $_R\mathcal{I}\subseteq ~_R\mathcal{F}_m$
by the assumption that $\sfli(R)\leq m$.
Conversely, let $_RM\in~_R\mathcal{F}_m\cap~{_R\mathcal{F}_m}{^{\perp}}$.
Consider the short exact sequence $$0 \to~ _RM \to ~_RI \to ~_RN \to 0$$ of $R$-modules
with $_RI$ injective.
Since $_RI$ and $_RM$ have finite flat dimension, so has $_RN$.
The assumption $\mathrm{FFD}(R)\leq m$ yields that $_RN$ is in $_R\mathcal{F}_m$, and so one has $\Ext_R^{1}(_RN,~_RM)=0$.
Thus the above exact sequence is split, and hence $_RM$ is injective.

(3)$\Rightarrow$(2). Suppose that $(_R\mathcal{F}_m~, ~{_R\mathcal{F}_m}^{\perp})$ is an injective cotorsion pair in $R\text{-}\Mod$.
Then ${_R\mathcal{F}_m}^{\perp} \subseteq ~_R\mathcal{GI}$ by \cite[Theorem 5.2]{Gil2016}.
On the other hand,  \cite[Theorem 5.6]{SS2020} tells us that $(^{\perp}\mathcal{GI}, \mathcal{GI})$ is a complete cotorsion pair, it is injective by \cite[Theorem 5.2]{Gil2016}.
It follows from \cite[Corollary 3.3]{Gil2017} that $_R\mathcal{GI} \subseteq~ {_R\mathcal{F}_m}^{\perp}$.
So $(_R\mathcal{F}_m~, ~_R\mathcal{GI})$=$(_R\mathcal{F}_m~, ~{_R\mathcal{F}_m}^{\perp})$ is an injective cotorsion pair in $R\text{-}\Mod$.

(2)$\Rightarrow$(1). Assume that $(_R\mathcal{F}_m~, ~_R\mathcal{GI})$ is an injective cotorsion pair in $R\text{-}\Mod$.
To see (1), first let $_RI$ be any injective left $R$-module.
Using the completeness of the cotorsion pair $(_R\mathcal{F}_m~, ~_R\mathcal{GI})$,
one gets a short exact sequence of $R$-modules
$$0 \to ~_RG \to ~_RF \to ~_RI \to 0$$
with $_RF\in ~_R\mathcal{F}_m$ and $_RG\in_R\mathcal{GI}$.
The above exact sequence is split as $\Ext_R^1(_RI, ~_RG)=0$. Hence $\fd(_RI)\leq \fd(_RF)\leq m$.
It follows that $\sfli(R)\leq m$.

Now, let $M$ be an $R$-module with $\fd_R(M)<\infty$.
Then $\Ext_R^{1}(_RM, ~_RG)=0$ for any $G\in ~_R\mathcal{GI}$ by \cite[Corollary 3.3]{Gil2017} as $(_R\mathcal{F}_m~, ~_R\mathcal{GI})$ is an injective cotorsion pair in $R\text{-}\Mod$.
This yields that $_RM\in ~_R\mathcal{F}_m$, and so one has $\FFD(R) \leq m$.
\end{proof}

Given a left $R$-module $_RM$ (viewed it as a complex), according to \cite{C-Y2021}, its Gorenstein flat-cotorsion dimension, denoted by $\Gfcd(_RM)$, is defined as follows
$$ \Gfcd(_RM) = \sup \left \{ m\in \mathbb{Z}~
\begin{array}{ |l} _RW^{\bullet}~\text{is a dg-flat-cotorsion resolution}~\text{of}~_RM~\text{with}~\h^i(_RW^{\bullet})=0 \\
\text{for all}~i<m~\text{and with}~\mathrm{C}^{m}(_RW^{\bullet})~\text{Gorenstein flat-cotorsion}.
\end{array}
\right\}
$$
Here a \emph{dg-flat-cotorsion resolution} of a complex $_RX^{\bullet}$ in $\Ch(R\text{-}\Mod)$ is a dg-flat resolution of  $_RX^{\bullet}$ with cotorsion components.
And a left $R$-module $_RN$ is called \emph{Gorenstein flat-cotorsion} if there is a \emph{totally acyclic complex of flat-cotorsion}
left $R$-modules $_RG^{\bullet}$, that is, $_RG^{\bullet}$ is a complex of flat-cotorsion left $R$-modules with $\Hom_R(_RG^{\bullet}, ~_RH)$ and $\Hom_R(_RH, ~_RG^{\bullet})$ acyclic
for any flat-cotorsion left $R$-module $_RH$, such that $_RN=\Z^{0}(_RG^{\bullet})$.
According to \cite[Theorem 5.7]{C-Y2021}, we know that $\Gfcd(-)$ is a refinement of $\Gfd(-)$.

\begin{thm} \label{3-3}
Let $m$ be a nonnegative integer. Then the following are equivalent for any $R$:
\begin{enumerate}
  \item $R$ is a ring with $\Gwgldim(R)\leq m$.
  \item $R$ is a left virtually Gorenstein ring with $\sup\{\Gfcd(_RM)~\mid~_RM\in R\text{-}\Mod\}\leq m$.
  \item There is a complete and hereditary cotorsion triple $(_R\mathcal{GP}, ~_R\mathcal{F}_m~,~ _R\mathcal{GI})$.
  \item There is a complete and hereditary cotorsion pair $(_R\mathcal{GP}, ~_R\mathcal{F}_m)$.
  \item There is an equality $_R\mathcal{F}_m=~_R\mathcal{GP}^{\perp}$.
\end{enumerate}
\end{thm}
\begin{proof}(1)$\Rightarrow$(2) follows from \cite[Theorem A]{DLW2023} and \cite[Fact 3.1]{C-W2023}.

(2)$\Rightarrow$(5). Suppose that $R$ is a left virtually Gorenstein ring with $\sup\{\Gfcd(_RM)~\mid~_RM\in R\text{-}\Mod\}\leq m$.
By \cite[Theorem 3.2]{C-W2023}, the assumption $\sup\{\Gfcd(_RM)~\mid~_RM\in R\text{-}\Mod\}\leq m$ is equivalent to $\max\{\sfli(R), \mathrm{FFD}(R)\}\leq m$.
Hence, there is a complete and hereditary cotorsion pair $(_R\mathcal{F}_m~,~ _R\mathcal{GI})$ in $R\text{-}\Mod$ by Lemma \ref{3-3-0}.
Now, since $R$ is left virtually Gorenstein, we get that $$_R\mathcal{GP}^{\perp}=~{^{\perp}{_R\mathcal{GI}}}= ~ _R\mathcal{F}_m.$$

(1)$\Rightarrow$(3) holds by \cite[Lemmas 3.2 and 2.1]{WZ2023} and (3)$\Rightarrow$(4)$\Rightarrow$(5) are clear.

(5)$\Rightarrow$(1). Assume that $R$ is a ring with an equality $_R\mathcal{F}_m=~_R\mathcal{GP}^{\perp}$. Then according to \cite[Corollary 3.4(1)]{CS2021}, the pair $(_R\mathcal{GP}, ~_R\mathcal{F}_m)$ forms a cotorsion pair. In particular, the class $~_R\mathcal{F}_m $ is closed under direct products. It is certainly closed under direct limits and, by \cite[Lemma 9.1.4]{EJ2000}, it is also closed under pure submodules. Hence, the class $~_R\mathcal{F}_m $ is definable. Then, as $~_R\mathcal{F}\subseteq ~_R\mathcal{F}_m$, we follow that $$~\overline{_RR}=~\overline{_R\mathcal{F}}\subseteq ~\overline{_R\mathcal{F}_m} =~_R\mathcal{F}_m=~_R\mathcal{GP}^{\perp}.$$
Thus, \cite[Theorem 4.4 and Corollary 4.5]{SS2020} (or \cite[Lemma 3.3]{WYZ2023}) yields $_R\mathcal{GP}=~_R\mathcal{PGF}$, and therefore $(_R\mathcal{PGF}, ~_R\mathcal{F}_m)$ is a cotorsion pair.
It is complete (and hereditary) by \cite[Theorem 4.9]{SS2020}.
Now, for any left $R$-module $_RM$, there is a short exact sequence of left $R$-modules
$$0\to ~_RM\to ~_RF\to ~_RG\to 0$$
with $_RF\in ~_R\mathcal{F}_m$ and $_RG\in ~_R\mathcal{PGF}$.
Therefore, one concludes from \cite[Lemma 2.2]{WZ2023} that $\Gfd(_RM)\leq m$, and so $\Gwgldim(R)\leq m$.
\end{proof}

\begin{rem} {\rm
It is not known whether the invariants $\Gwgldim(R)$ and $\sup\{\Gfcd(_RM)~\mid~_RM\in R\text{-}\Mod\}$ always agree for an arbitrary ring $R$. By \cite[Section 5]{C-W2023} they coincide over a large class of rings. Theorem \ref{3-3} adds another class of rings for which the two invariants are the same: left virtually Gorenstein rings.}

\end{rem}

\begin{rem} \label{3-5}{\rm
Let $R$ be a  ring with $\Gwgldim(R)<\infty$. By Theorem \ref{3-3} and \cite[Lemma 2.1]{WZ2023} (see also \cite[Lemma 3.2]{WZ2023}), there are complete and hereditary cotorsion pairs
$(_R\mathcal{GP},~_R\mathcal{F}^\wedge)$ and $(_R\mathcal{F}^\wedge, ~_R\mathcal{GI})$  in $R$-Mod.
By using \cite[Theorems 5.2 and 5.4]{Gil2016}, such cotorsion pairs are further projective and injective respectively.
Now the induced cotorsion pairs $(\dg(_R\mathcal{GP}), \widetilde{_R\mathcal{F}^\wedge})$ and $(\widetilde{_R\mathcal{F}^\wedge},~\dg(_R\mathcal{GI}))$ in $\Ch(R\text{-}\Mod)$
are projective and injective respectively by \cite[Propositions 7.2 and 7.3]{Gil2016}.}
\end{rem}

Let $R$ be any ring. As already mentioned, by \cite[Theorem 5.6]{SS2020}, there is a complete and hereditary cotorsion pair $(^{\perp}{_R\mathcal{GI}},~_R\mathcal{GI})$ in $R$-Mod.
So we have the notion of dg-Gorenstein injective complex (see Definition \ref{10-9}). As usual, we denote by $_R\mathcal{C}ot$
 the subcategory of $R$-Mod consisting of all \emph{cotorsion} left $R$-modules $_RM$
(that is, $_RM\in ~_R\mathcal{F}^{\perp}$); therefore, Ch$(_R\mathcal{C}ot)$ denotes the subcategory of \Ch$(R\text{-}\Mod)$ consisting of all cotorsion modules in our sense. The next result shows that dg-Gorenstein injective complexes are particularly simple over rings of finite Gorenstein weak global dimension.

\begin{lem} \label{4}Let $R$ be a a ring with $\Gwgldim(R)<\infty$. Then  the following are equivalent for any $_RM^{\bullet}\in \Ch(R\text{-}\Mod)$:

\begin{enumerate}
  \item $_RM^{\bullet}$ is dg-Gorenstein injective.
  \item $_RM^{\bullet}\in \Ch(R\text{-}\GInj)$.
  \item Each $_RM^{i}$ is Gorenstein injective.
\end{enumerate}

 \end{lem}

\begin{proof} (2)$\Leftrightarrow$(3) follows from \cite[Proposition 2.8]{YL2011} and (1)$\Rightarrow$(2) is clear.

(2)$\Rightarrow$(1). Suppose that $R$ is a ring with $\Gwgldim(R)=m<\infty$.
Let $_RM^{\bullet}\in \Ch(R\text{-}\GInj)$. Then $_RM^{\bullet}\in \Ch(_R\mathcal{C}ot)$ by \cite{Sto2014}.
Due to Remark \ref{3-5} there exists a short exact sequence in $\Ch(R\text{-}\Mod)$
$$0\to ~_RM^{\bullet}\to ~_RT^{\bullet}\to ~_RN^{\bullet}\to0$$
with $_RT^{\bullet}\in \dg(_R\mathcal{GI})$ and $_RN^{\bullet}\in \widetilde{_R\mathcal{F}_m}$.
Note that $_RN^{\bullet}$ is in $\Ch(_R\mathcal{C}ot)$ since $_RM^{\bullet}$ and $_RT^{\bullet}$ are so. Now $_RN^{\bullet} \in \widetilde{_R\mathcal{F}_m}\cap \Ch(_R\mathcal{C}ot)$ is contractible by \cite[Theorem 5.2(2)]{BCE2020}. Thus, $_RM^{\bullet}$ and $_RT^{\bullet}$ are homotopically equivalent, and hence $_RM^{\bullet}$ is dg-Gorenstein injective.
\end{proof}

As mentioned in Remark \ref{3-5}, there is a complete and hereditary cotorsion pair $(_R\mathcal{GP},~_R\mathcal{F}^\wedge)$ whenever $\Gwgldim(R)<\infty$. Now one has the notion of dg-Gorenstein projective complexes (see Definition \ref{10-9}) in this case. As in the previous case, dg-Gorenstein projective complexes can be characterized as follows:

\begin{lem} \label{5}Let $R$ be a a ring with $\Gwgldim(R)<\infty$. Then  the following are equivalent for any $_RM^{\bullet}\in \Ch(R\text{-}\Mod)$:

\begin{enumerate}
  \item $_RM^{\bullet}$ is dg-Gorenstein projective.
  \item $_RM^{\bullet}\in \Ch(R\text{-}\GProj)$.
  \item Each $_RM^{i}$ is Gorenstein projective.
\end{enumerate}

 \end{lem}

\begin{proof}

(2)$\Leftrightarrow$(3) follows from \cite[Theorem 2.2]{YL2011} and (1)$\Rightarrow$(2) is obvious.

(2)$\Rightarrow$(1). Assume that $R$ is a ring with $\Gwgldim(R)=m<\infty$.
Let $_RM^{\bullet}\in \Ch(R\text{-}\GProj)$.
By Remark \ref{3-5} there exists a short exact sequence in $\Ch(R\text{-}\Mod)$
$$0\to ~_RK^{\bullet}\to ~_RT^{\bullet}\to ~_RM^{\bullet}\to0$$
with $_RT^{\bullet}\in \dg(_R\mathcal{GP})$ and $_RK^{\bullet}\in \widetilde{_R\mathcal{F}_m}$.
Since both $_RM^{\bullet}$ and $_RT^{\bullet}$ belong to $\Ch(R\text{-}\GProj)$, one gets $_RK^{\bullet}\in\Ch(R\text{-}\GProj)$, and hence $_RK^{\bullet} \in \Ch(R\text{-}\GProj)\cap \widetilde{_R\mathcal{F}_m}$.
From Theorem \ref{3-3} we have that $_R\mathcal{F}_m = ~_R\mathcal{GP}^\perp$. This implies that each $K^j\in ~_R\mathcal{GP}\cap ~_R\mathcal{GP}^\perp=~_R\mathcal{P} $. Therefore, $_RK^{\bullet}$ is an acyclic complex in $\Ch(R\text{-}\Proj)$  with cycle modules in $_R\mathcal{F}_m$.
Then, it is easy to check that $_RK^{\bullet}$ is an acyclic complex in $\Ch(R\text{-}\Proj)$  with cycle modules in $_R\mathcal{F}$. Now $_RK^{\bullet}$ is contractible by a result \cite[Theorem 2.5]{BG2000} due to Benson and Goodearl (see also \cite[Corollary 3.6(2)]{BCE2020}). Thus,  $_RM^{\bullet}$ and $_RT^{\bullet}$ are homotopically equivalent, and so $_RM^{\bullet}$ is dg-Gorenstein projective.\end{proof}

Now, we are in a position to state and prove our main result. For a ring $R$, recall that $\K(\Proj\text{-}R)$ (resp. $\K(\Inj\text{-}R)$, $\K(\GProj\text{-}R)$ and $\K(\GInj\text{-}R)$) denotes
the homotopy category of all projective (resp. injective, Gorenstein projective and Gorenstein injective) right $R$-modules,
and that their homotopy subcategory consisting of all exact complexes will be denoted by $\K_{\mathrm{ac}}(\Proj\text{-}R)$ , $\K_{\mathrm{ac}}(\Inj\text{-}R)$, $\K_{\mathrm{ac}}(\GProj\text{-}R)$ and $\K_{\mathrm{ac}}(\GInj\text{-}R)$
respectively.

\begin{thm} \label{8}
Let $R$ be  a ring with $\Gwgldim(R)<\infty$.
 Then there is a triangle-equivalence
$\K(R\text{-}\GProj)$ $\simeq \K(R\text{-}\GInj)$
which restricts to a triangle-equivalence $\K(R\text{-}\Proj)\simeq \K(R\text{-}\Inj).$
 \end{thm}
\begin{proof}
Suppose that $R$ is a ring with $\Gwgldim(R)<\infty$. Then it is a left (and right) virtually Gorenstein  ring with a complete and hereditary triple $(_R\mathcal{GP}, ~_R\mathcal{F}^\wedge, ~_R\mathcal{GI})$ by Theorem \ref{3-3}.
So Example \ref{T-9}(2) tells us the following triangle-equivalences
$$\K(\dg(_R\mathcal{GP}))\simeq\mathbf{D}_{\mathrm{GP}}(R\text{-}\Mod)\simeq\mathbf{D}_{\mathrm{GI}}(R\text{-}\Mod)\simeq\K(\dg(_R\mathcal{GI})).$$
Using Lemmas \ref{4} and \ref{5}, there are equalities
$$\K(\dg(_R\mathcal{GP}))=\K(R\text{-}\GProj)\quad\text{and}\quad\K(\dg(_R\mathcal{GI}))=\K(R\text{-}\GInj).$$
Thus, we have a triangle-equivalence $\K(R\text{-}\GProj)\simeq \K(R\text{-}\GInj)$.
Furthermore, if we denote by $F:\K(R\text{-}\GProj)\to \K(R\text{-}\GInj)=$ $\K(\dg(_R\mathcal{GP}))\to \K(\dg(_R\mathcal{GI}))$ and $F':\K(R\text{-}\GInj)\to \K(R\text{-}\GProj)=$ $\K(\dg(_R\mathcal{GI}))\to \K(\dg(_R\mathcal{GP}))$ the associate equivalent functor respectively,
then the rules of $F$ and $F'$ on objects are defined as follows (see Example \ref{T-9}(2)):
\begin{center}$\begin{aligned}&F: ~_RX^{\bullet}\mapsto ~_RY^{\bullet}, \text{where}~_RY^{\bullet}~\text{is a dg-Gorenstein injective resolution of}~_RX^{\bullet}; \\
&F': ~_RN^{\bullet}\mapsto ~_RM^{\bullet}, \text{where}~_RM^{\bullet}~\text{is a dg-Gorenstein projective resolution of}~_RN^{\bullet}.
\end{aligned}$\end{center}
To see the restricted triangle-equivalence, it suffices to show the essential image of $F$ (resp. $F'$) is contained in $\K(R\text{-}\Inj)$ (resp. $\K(R\text{-}\Proj)$).
For this, let $_RX^{\bullet}$ (resp. $_RN^{\bullet}$) be any complex in $\K(R\text{-}\Proj)$ (resp. $\K(R\text{-}\Inj)$).
Then by Remark \ref{3-5}, there are short exact sequences in $\Ch(R\text{-}\Mod)$
$$0\to ~_RX^{\bullet} \overset{\alpha^{\bullet}}\to ~_RG^{\bullet}\to ~_RC^{\bullet}\to 0\quad\text{and}\quad 0\to ~_RK^{\bullet}\to ~_RH^{\bullet}\overset{\beta^{\bullet}}\to ~_RN^{\bullet}\to 0$$
 with $_RG^{\bullet}\in \dg(_R\mathcal{GI})$, $_RH^{\bullet}\in \dg(_R\mathcal{GP})$ and $_RC^{\bullet}, ~_RK^{\bullet}\in~\widetilde{_R\mathcal{F}^\wedge}$.
Note that $\alpha^{\bullet}$ and $\beta^{\bullet}$ are quasi-isomorphisms as $_RC^{\bullet}$ and $_RK^{\bullet}$ are exact complexes, and hence, $_RG^{\bullet}$ is a dg-Gorenstein injective resolution of $_RX^{\bullet}$
and $~_RH^{\bullet}$ is a dg-Gorenstein projective resolution of $_RN^{\bullet}$.
Since each $_RC^{k}$ is in $^{\perp}{_R\mathcal{GI}}$ and each $_RX^{k}$ is projective (of course in $^{\perp}{_R{\mathcal{GI}}}$),
and that each $_RK^{k}$ is in $_R\mathcal{GP}^{\perp}$ and each $_RN^{k}$ is injective (of course in $_R\mathcal{GP}^{\perp}$),
we conclude that $_RG^{k}\in ~_R\mathcal{GI}\cap~{^{\perp}{_R\mathcal{GI}}}=~_R\mathcal{I}$ for all $k\in \mathbb{Z}$
and that $_RH^{k}\in {_R\mathcal{GP}}\cap~{_R\mathcal{GP}^{\perp}}=~_R\mathcal{P}$ for all $k\in \mathbb{Z}$.
Therefore, $F(_RX^{\bullet})=~_RG^{\bullet}\in \K(R\text{-}\Inj)$ and $F'(_RN^{\bullet})=~_RH^{\bullet}\in \K(R\text{-}\Proj)$, as desired.
\end{proof}

\begin{rem} \label{9}{\rm
Let $R$ be a ring with $\Gwgldim(R)<\infty$.
Applying analogous arguments to those in the proof of Theorem \ref{8}, we can see that the triangle-equivalence
$\K(R\text{-}\GProj)\simeq \K(R\text{-}\GInj)$ restricts to a triangle-equivalence $\K_{\mathrm{ac}}(R\text{-}\Proj)\simeq \K_{\mathrm{ac}}(R\text{-}\Inj)$
(resp. $\K_{\mathrm{ac}}(R\text{-}\GProj)\simeq \K_{\mathrm{ac}}(R\text{-}\GInj)$).}
 \end{rem}

Let us consider two immediate consequences of Theorem \ref{8}. The first one is about left Gorenstein rings. In particular, it recovers \cite[Theorem B]{Chen2010}.

\begin{cor} \label{9-0-0}
Let $R$ be a left Gorenstein ring. The following assertions hold:
\begin{enumerate}
\item There is a triangle-equivalence
$\K(R\text{-}\GProj)$ $\simeq \K(R\text{-}\GInj)$
which restricts to triangle-equivalences  $\K(R\text{-}\Proj)\simeq \K(R\text{-}\Inj)$, $\K_{\mathrm{ac}}(R\text{-}\Proj)\simeq \K_{\mathrm{ac}}(R\text{-}\Inj)$ and $\K_{\mathrm{ac}}(R\text{-}\GProj)\simeq \K_{\mathrm{ac}}(R\text{-}\GInj)$.
\item There is a triangle-equivalence
$\K(\GProj\text{-}R)$ $\simeq \K(\GInj\text{-}R)$
which restricts to triangle-equivalences  $\K(\Proj\text{-}R)\simeq \K(\Inj\text{-}R)$, $\K_{\mathrm{ac}}(\Proj\text{-}R)\simeq \K_{\mathrm{ac}}(\Inj\text{-}R)$ and $\K_{\mathrm{ac}}(\GProj\text{-}R)\simeq \K_{\mathrm{ac}}(\GInj\text{-}R)$
\end{enumerate}

\end{cor}

\begin{proof} Due to \cite[Theorem 3.7]{WYZ2023} or \cite[Theorem 5.13]{C-W2023}, we know that $\Gwgldim(R)\leq \Ggldim(R)$.
Thus, any left Gorenstein ring $R$ (i.e., $R$ is of  $\Ggldim(R)<\infty$) admits  $\Gwgldim(R)<\infty$.
Now, the first statement follows from Theorem \ref{8} and Remark \ref{9}. For the second statement, simply note that Gorenstein weak global dimension is left-right symmetric (see \cite{CET2021}).
\end{proof}

In what follows, we always let $n\in \mathbb{N}\cup\{\infty\}$ and denote by $_R\mathcal{FP}_n\mathcal{F}$ (resp $_R\mathcal{FP}_n\mathcal{I}$) the subcategory of $R$-Mod
consisting of all FP$_n$-flat (resp. FP$_n$-injective) modules.
Let us recall some notions from \cite{BP2017}.

(1) A left $R$-module~$_RF$ is called \emph{finitely} $n$-\emph{presented}, if~$_RF$ has a projective resolution of the form
\begin{center}$ _RP_{n}\rightarrow
\cdots \rightarrow
~_RP_{1}\rightarrow ~_RP_{0}\rightarrow~ _RM\rightarrow0$ \end{center}
in which each $_RP_i$ is finitely generated projective.

(2) A left $R$-module $_RM$ (resp. a right $R$-module $N_R$) is said to be \emph{FP}$_n$-\emph{injective} (resp. \emph{FP}$_n$-\emph{flat})
if $\Ext_R^{1}(_RM,~_RF)=0$ (resp. $\Tor^{R}_1(N_R,~_RF)=0$) for all finitely $n$-presented left $R$-modules $_RF$.

(3) A ring $R$ is called left (resp. right) $n$-\emph{coherent} if each finitely $n$-presented left (resp. right) $R$-module is finitely $n+1$-presented.

Note that the notion of FP$_n$-injective module unifies the notions of injective, FP-injective and \emph{absolutely clean} \cite{BGH2014} modules, in the sense of that
FP$_0$-injective (resp. FP$_1$-injective and FP$_\infty$-injective)-modules are exactly injective (resp. FP-injective and absolutely clean) modules.
Dually,  the notion of FP$_n$-flat modules unifies the notions of flat and \emph{level} \cite{BGH2014} modules, in the sense of that
FP$_0$-flat or FP$_1$-flat (resp. and FP$_\infty$-flat)-modules are exactly flat (resp. level) modules.
On the other hand, the notion of left $n$-coherent ring unifies the notions of left Noether and left coherent rings, in the sense of that
left 0-coherent (resp. left 1-coherent) rings are exactly left Noether (resp. left coherent) rings. Meanwhile, arbitrary rings are left (right) $\infty$-coherent rings.

In order to give a negative answer to a question posed by Gillespie in \cite{Gil2017}, Wang, Liu and Yang \cite{WLY2018} introduced the following notion.

\begin{df} \label{9-1}
{\rm (=\cite[Theorem 4.1]{WLY2018}). Let $n\in \mathbb{N}\cup\{\infty\}$.
A ring $R$ is called \emph{Gorenstein} $n$-\emph{coherent} if $R$ is a two-sided $n$-coherent ring with the condition $(m\text{-}\mathrm{FP}_n\mathrm{I}\text{-}\mathrm{F})$
for some nonnegative integer $m$. Here, for any integer $m\geq0$,
$R$ is of \emph{the condition} \emph{(m\text{-}FP}$_n$\emph{I\text{-}F)} if all left and right FP$_n$-injective $R$-modules have finite flat dimension at most $m$.
}
\end{df}

Recall that a ring $R$ is \emph{Gorenstein}~\cite{Iwa1979}
(resp. \emph{Ding-Chen}~\cite{DC1996, Gil2010}) if $R$ is an $m$-Gorenstein ring (resp. $m$-FC ring) for some nonnegative integer $m$, i.e., $R$
is a two-sided noetherian (resp. two-sided coherent) ring with self-injective (resp. self-FP-injective)
dimension at most $m$ on both sides.

Just as the notion of left (right) $n$-coherent rings unifies the notions of left (right) noetherian rings and left (right) coherent rings,
the notion of Gorenstein $n$-coherent rings unifies the notions of Gorenstein rings and Ding-Chen rings, in the sense of that
Gorenstein 0-coherent (resp. Gorenstein 1-coherent) rings are exactly Gorenstein (resp. Ding-Chen) rings.



Now we can state the second consequence of Theorem \ref{8}.

\begin{cor} \label{9-2}
Let $n\in \mathbb{N}\cup\{\infty\}$ and $R$ be a Gorenstein $n$-coherent ring.
 Then there is a triangle equivalence
$\K(R\text{-}\GProj)\simeq \K(R\text{-}\GInj)$
which
 restricts to triangle equivalences $$\K(R\text{-}\Proj)\simeq \K(R\text{-}\Inj)\quad\text{and}\quad\K_{\mathrm{ac}}(R\text{-}\Proj)\simeq \K_{\mathrm{ac}}(R\text{-}\Inj).$$
 \end{cor}
\begin{proof} Let $R$ be a Gorenstein $n$-coherent ring $R$ with the condition ($m$ $\text{-FP}_n$I$\text{-}$F) for $m<\infty$.
Then all FP$_n$-injective left and right $R$-modules  have finite flat dimension at most $m$.
In particular, all injective left and right $R$-modules have finite flat dimension at most $m$.
Now \cite[Theorem 5.3]{Emm2012} shows that $R$ admits $\Gwgldim(R)\leq m$, and then Theorem \ref{8} and Remark \ref{9} apply.
\end{proof}

The particular case of $n=0,1$ in Corollary \ref{9-2} is as follows:

\begin{cor} \label{9-2-1}
Let $R$ be a Gorenstein ring or more generally a Ding-Chen ring. Then the triangle-equivalence $\K(R\text{-}\GProj)\simeq \K(R\text{-}\GInj)$ holds true, moreover, it
 restricts to triangle equivalences $$\K(R\text{-}\Proj)\simeq \K(R\text{-}\Inj)\quad\text{and}\quad\K_{\mathrm{ac}}(R\text{-}\Proj)\simeq \K_{\mathrm{ac}}(R\text{-}\Inj).$$
 \end{cor}

We end this section by some examples, which show that there are rings with finite Gorenstein weak global dimension which are not left Gorenstein.

\begin{exa} \label{9-6}{\rm
Let $R=F_{\alpha}$ be the free Boolean ring on $\aleph_{\alpha}$ generators with $\alpha$ an infinite cardinality.
By \cite[Example 3.3]{Wz2017}, $R$ is a commutative von Neumann regular ring (hence a
commutative Ding-Chen ring with $\Gwgldim(R) = \wgldim(R) = 0$) that has infinite
Gorenstein global dimension.

Let $x_1, x_2, \cdots, x_m$ be variables which commute with each element of $R$.
We consider the the polynomial rings in such variables with coefficients in $R$.
By \cite[Theorems 4.8 and 4.13]{LYZ2023}, one has
$$\Gwgldim(R[x_i])=\Gwgldim(R)+1=1 \quad\text{and}\quad \Ggldim(R[x_i])=\Ggldim(R)+1=\infty.$$
for each $1\leq i\leq m$.
Now by induction, one gets that
$\Gwgldim(R[x_1, \cdots, x_k])=\Gwgldim(R)+k=k$ and $\Ggldim(R[x_1,\cdots, x_k])=\Ggldim(R)+k=\infty.$
for each $1\leq k\leq m$.

Now,  both $R$ and each $R[x_i]$ as well as each $R[x_1, \cdots, x_k]$
admit finite Gorenstein weak global dimension. However, none of them are left Gorenstein.
}
 \end{exa}

\section{\bf Homotopy equivalence on totally acyclic complexes}

Let $R$ be a ring with finite Gorenstein weak global dimension.
This section is devoted to give a link between $\K_{\mathrm{tac}}(R\text{-}\Proj)$
and $\K_{\mathrm{tac}}(R\text{-}\Inj)$, the homotopy categories of totally acyclic complexes of
projective left $R$-modules and totally acyclic complexes of injective left $R$-modules.

\begin{lem} \label{9-6-1}  Let $R$ be a ring with $\Gwgldim(R)<\infty$. Then there are equalities
$$\K_{\mathrm{tac}}(R\text{-}\Proj)=\K_{\mathrm{ac}}(R\text{-}\Proj)\quad\text{and}\quad\K_{\mathrm{tac}}(R\text{-}\Inj)=\K_{\mathrm{ac}}(R\text{-}\Inj).$$
\end{lem}

\begin{proof}
Using  \cite[Theorems 2.3 and 2.6]{WZ2023}, we get that  any exact complex of projective (injective) left $R$-modules is totally acyclic.
In other words, the containments $$\K_{\mathrm{tac}}(R\text{-}\Proj)\supseteq \K_{\mathrm{ac}}(R\text{-}\Proj)\quad\text{and}\quad\K_{\mathrm{tac}}(R\text{-}\Inj)\supseteq \K_{\mathrm{ac}}(R\text{-}\Inj)$$
hold true. On the other hand, the converse containments are trivial.
\end{proof}

Now we can give a link between $\K_{\mathrm{tac}}(R\text{-}\Proj)$
and $\K_{\mathrm{tac}}(R\text{-}\Inj)$.

\begin{prop} \label{9-6-2}  Let $R$ be a ring with $\Gwgldim(R)<\infty$. Then there are triangle equivalences
$$\K_{\mathrm{tac}}(R\text{-}\Proj)=\K_{\mathrm{ac}}(R\text{-}\Proj)\simeq \K_{\mathrm{ac}}(R\text{-}\Inj)=\K_{\mathrm{tac}}(R\text{-}\Inj).$$
\end{prop}

\begin{proof} The equalities follows from Lemma \ref{9-6-1} and the triangle equivalence from Remark \ref{9}.
\end{proof}

It is well known that over an arbitrary ring $R$ the subcategory $_R\mathcal{GP}$ (resp., $_R\mathcal{GI}$)  is a Frobenius category with projective-injective objects projective (resp. injective) left $R$-modules.
Hence, the stable categories $\underline{_R\mathcal{GP}}$ and $\underline{_R\mathcal{GI}}$ are triangulated categories.

\begin{cor} \label{9-6-3}  Let $R$ be a ring with $\Gwgldim(R)<\infty$. Then there is a triangle equivalence
$$\underline{_R\mathcal{GP}}\simeq \underline{_R\mathcal{GI}}.$$
\end{cor}

\begin{proof} According to \cite[Proposition 7.2 and its dual]{Kra2005}, we know that there are triangle equivalences
$$\K_{\mathrm{tac}}(R\text{-}\Proj)\simeq \underline{_R\mathcal{GP}}\quad\text{and}\quad\K_{\mathrm{tac}}(R\text{-}\Inj)\simeq \underline{_R\mathcal{GI}}.$$
Now Proposition \ref{9-6-2} applies.
\end{proof}

We end this section with the next corollary.

\begin{cor} \label{9-6-4}  Let $R$ be a left Gorenstein ring or a Gorenstein $n$-coherent ring for some $n\in \mathbb{Z}\cup\{\infty\}$. Then there are  triangle equivalences
$$\underline{_R\mathcal{GP}}\simeq\K_{\mathrm{tac}}(R\text{-}\Proj)=\K_{\mathrm{ac}}(R\text{-}\Proj)\simeq \K_{\mathrm{ac}}(R\text{-}\Inj)=\K_{\mathrm{tac}}(R\text{-}\Inj)\simeq \underline{_R\mathcal{GI}}.$$ and $$\underline{\mathcal{GP}_R}\simeq\K_{\mathrm{tac}}(\Proj\text{-} R)=\K_{\mathrm{ac}}(\Proj \text{-}R)\simeq \K_{\mathrm{ac}}(\Inj\text{-}R )=\K_{\mathrm{tac}}(\Inj\text{-}R )\simeq \underline{\mathcal{GI}_R}.$$
\end{cor}

\begin{proof} By the proofs of Corollaries \ref{9-0-0} and \ref{9-2}, we know that $\Gwgldim(R)<\infty$.
Now the triangle equivalences hold by Proposition \ref{9-6-2} and \cite[Proposition 7.2 and its dual]{Kra2005}.
\end{proof}

\section{\bf Homotopy equivalence with respect to duality pairs}

Let $R$ be a flat-typed $(\mathcal{L},\mathcal{A})$-Gorenstein ring with respect to a bi-complete
duality pair $(\mathcal{L},\mathcal{A})$ (see Definitions 2.1-2.3 for details). In this section,  we will obtain a homotopic equivalence which relates the subcategory $\mathcal{A}$ and the subcategory $\mathcal{F}$ of flat modules (see Theorem \ref{7-5}).
Moreover, we will prove that such equivalent homotopy subcategories are compactly generated under some mild condition (see Theorem \ref{7-7}).

For a ring with a bi-complete duality pair $(\mathcal{L},\mathcal{A})$, we denote by $\Ch(_R\mathcal{A}\text{-}\Inj)$~(resp. $\Ch(_R\mathcal{L}\text{-}\Proj)$) the subcategory of $\Ch(R\text{-}\Inj)$~(resp. $\Ch(R\text{-}\Proj)$) consisting of all complexes $_RI^{\bullet}$ such that any morphism $_RA^{\bullet}\to ~_RI^{\bullet}$ is null-homotopic, whenever $_RA^{\bullet} \in \widetilde{_R\mathcal{A}}$  (resp. $_RP^{\bullet}$ such that any morphism $_RP^{\bullet}\to ~_RL^{\bullet}$ is null-homotopic, whenever $_RL^{\bullet} \in \widetilde{_R\mathcal{L}}$). And then denote by $\K(_R\mathcal{A}$-\Inj$)$ and $\K(_R\mathcal{L}$-\Proj$)$ the corresponding homotopy subcategories, respectively.

\begin{lem} \label{7-0-0} Let $R$ be a ring. Then we have a triangle equivalence $$\K(_R\mathcal{F}\cap ~_R\mathcal{C}ot)=\K(\Ch(_R\mathcal{F})\cap \dg(_R\mathcal{C}ot))\simeq \K(R\text{-}\Proj).$$
\end{lem}
\begin{proof} Note from \cite[Theorem 5.3]{BCE2020} that dg$(_R\mathcal{C}ot)=$ Ch$(_R\mathcal{C}ot)$. Now the equality holds.
The remaining triangle equivalence comes from \cite[Corollary 7.6]{Gil16}.
\end{proof}

Now we can consider some triangle-equivalences involving  $\K(_R\mathcal{L}\text{-}\Proj)$.

\begin{prop} \label{7-4}Let $R$ be a flat-typed $(\mathcal{L},\mathcal{A})$-Gorenstein ring with respect to a bi-complete
duality pair $(\mathcal{L},\mathcal{A})$.
Then we have triangle-equivalences
$$\K(_R\mathcal{L}\text{-}\Proj)=\K(R\text{-}\Proj)\simeq \K(_R\mathcal{F}\cap ~_R\mathcal{C}ot).$$
 \end{prop}
\begin{proof}
The triangle-equivalence is obtained by Lemma \ref{7-0-0}.
To see the equality, it is obvious that $\K(_R\mathcal{L}\text{-}\Proj)\subseteq \K(R\text{-}\Proj)$.
Conversely, let $_RX^{\bullet}\in \K(R\text{-}\Proj)$.
For any complex $_RL^{\bullet}\in \widetilde{_R\mathcal{L}}$, one has $_RL^{\bullet}$ is an exact complex with each cycle module $Z^{k}(_RL^{\bullet})$ in $_R\mathcal{L}$.
Note from \cite[Proposition 4.5]{WD2020} that any module in $_R\mathcal{L}$ has finite flat dimension in $R$-Mod.
Hence $_RL^{\bullet}$ has finite flat dimension in \Ch$(R\text{-Mod})$, and so, there exists an integer $m\geq0$ and  an exact sequence in \Ch$(R\text{-Mod})$
$$0\to ~_RF_m^{\bullet}\to \cdots\to ~_RF_1^{\bullet}\to ~_RF_0^{\bullet}\to ~_RL^{\bullet}\to0$$
 with each $_RF_i^{\bullet}$ a flat complex. This induces the following exact sequence in \Ch$(\mathbb{Z}\text{-Mod})$
$$0\to \mathcal{H}om(_RP^{\bullet},~_RF_m^{\bullet})\to \cdots\to \mathcal{H}om(_RP^{\bullet},~_RF_1^{\bullet})\to \mathcal{H}om(_RP^{\bullet},~_RF_0^{\bullet})\to \mathcal{H}om(_RP^{\bullet},~_RL^{\bullet})\to 0$$
for any $_RP^{\bullet}\in \Ch(R\text{-}\Proj)$.
Note from \cite[Theorem 8.6]{Nee2008} that each morphism $f_{i}^{\bullet}:~_RP^{\bullet}\to ~~_RF_i^{\bullet}$ is null-homotopic,
i.e., each $\mathcal{H}om(_RP^{\bullet},~_RF_i^{\bullet})$ is exact. And then so is $\mathcal{H}om(_RP^{\bullet}, ~_RL^{\bullet})$.
This in turn shows that any morphism $f^{\bullet}:~_RP^{\bullet}\to ~_RL^{\bullet}$ with $_RL^{\bullet}\in \widetilde{_R\mathcal{L}}$ is null-homotopic.
So, by taking $~_RP^{\bullet}= ~_RX^{\bullet}$ in the above, we conclude that $_RX^{\bullet}\in\K(_R\mathcal{L}\text{-}\Proj)$.
\end{proof}

Given a left $R$-module $_RM$ and an integer $k\in \mathbb{Z}$, following \cite{Gil201702}, we denote by D$^k(_RM)$ the complex $\cdots \to 0 \to~_RM \overset{Id}
\to~_RM \to0\to\cdots $ with $_RM$ in the $k$ and $k+1$th degrees.

For a bi-complete duality pair $(\mathcal{L},\mathcal{A})$, there are two exact subcategories $_R\mathcal{L}$ and $_R\mathcal{A}$ of $R$-Mod.
They induce two exact subcategories $\Ch(_R\mathcal{L})$ and $\Ch(_R\mathcal{A})$ of \Ch$(R\text{-Mod})$, with their inherited degreewise exact structures respectively.

Note that the projective and injective cotorsion pairs can be also defined in a weakly idempotent complete exact category (see \cite{Gil16}). Furthermore, the results \cite[Propositions 3.6 and 3.7] {Gil2016} for the characterization of such cotorsion pair still hold in this setting.

\begin{lem} \label{7} Let $R$ be a ring with a bi-complete
duality pair $(\mathcal{L},\mathcal{A})$.
Then we have a triangle-equivalence $\K(_R\mathcal{A}\text{-}\Inj)\simeq \K(\dg(^{\perp}{_R\mathcal{A}})\cap \Ch(_R\mathcal{A}))$.
 \end{lem}
\begin{proof} By \cite[Corollary 4.5]{Gil16}, it suffices to show that the triple $(\K(\dg(^{\perp}{{_R\mathcal{A}}})\cap \Ch(_R\mathcal{A})), ~\widetilde{_R\mathcal{A}}, \K(R\text{-}\Inj))$
is a localizing cotorsion triple (in the sense of \cite{Gil16}) in $\Ch(_R\mathcal{A})$.
It suffices to see that there is a cotorsion triple $$(\dg(^{\perp}{_R\mathcal{A}})\cap \Ch(_R\mathcal{A}),~\widetilde{_R\mathcal{A}},~\Ch(_R\mathcal{A}\text{-}\Inj))$$
in the exact category $\Ch(_R\mathcal{A})$ such that the left-hand cotorsion pair is projective  and the right-hand one is injective.
For this, we prove that

{\bf Fact 1}: There is an injective cotorsion pair $(_R\mathcal{W}^{\bullet},~\Ch(_R\mathcal{A}\text{-}\Inj))$ in \Ch$(R\text{-Mod})$, where $_R\mathcal{W}^{\bullet}$ stands for $^{\perp}{\Ch(_R\mathcal{A}\text{-}\Inj)}$.

Note that the subcategory $\widetilde{_R\mathcal{A}}$ is closed under pure subcomplexes and pure quotients of complexes
(since one can check that the pair $(\widetilde{\mathcal{L}},\widetilde{\mathcal{A}})$ is a bi-complete duality in \Ch$(R\text{-Mod})$, see \cite{EIP2020}).
Hence, there is a set $\{_RA^{\bullet}_{\lambda}~|~\lambda\in \Lambda\}$ of complexes in $\widetilde{_R\mathcal{A}}$ such that any complex in $\widetilde{_R\mathcal{A}}$  is a transfinite extension of complexes in  $\{_RA^{\bullet}_{\lambda}~|~\lambda\in \Lambda\}$.
Now we set $$_R\mathcal{S}^{\bullet}=\{\D^k(_RR/_RI)~|~k\in \mathbb{Z},~_RI~\text{is a left ideal of}~R\}\cup \{_RA^{\bullet}_{\lambda}~|~\lambda\in \Lambda\}.$$
One can verify  that ${_R\mathcal{S}^{\bullet}}^{\perp}=\Ch(_R\mathcal{A}\text{-}\Inj)$.
Therefore, the pair $(_R\mathcal{W}^{\bullet},~\Ch(_R\mathcal{A}\text{-}\Inj))$ is a complete cotorsion pair in \Ch$(R\text{-Mod})$.
In addition, $_R\mathcal{W}^{\bullet}$ is closed under direct summands in \Ch$(R\text{-Mod})$.
Let $$(\dag)\quad0\to ~_RX^{\bullet}\to~ _RY^{\bullet}\to ~_RZ^{\bullet}\to0$$ be a short exact sequence in \Ch$(R\text{-Mod})$ with $_RX^{\bullet}, ~_RY^{\bullet}\in ~_R\mathcal{W}^{\bullet}$.
For any $_RA^{\bullet}\in \Ch(_R\mathcal{A}\text{-}\Inj)\subseteq \Ch(R\text{-}\Inj)$, one has a short exact sequence in \Ch$(\mathbb{Z}\text{-Mod})$
$$0\to \mathcal{H}om(_RZ^{\bullet},~_RA^{\bullet})\to \mathcal{H}om(_RY^{\bullet},~_RA^{\bullet})\to \mathcal{H}om(_RX^{\bullet},~_RA^{\bullet})\to 0.$$
Note that $\mathcal{H}om(_RZ^{\bullet},~_RA^{\bullet})$ is exact since the assumption on $_RX^{\bullet}$ and $_RY^{\bullet}$
yields that both $\mathcal{H}om(_RX^{\bullet},~_RA^{\bullet})$ and $\mathcal{H}om(_RY^{\bullet},~_RA^{\bullet})$ are exact.
It follows that $\Ext_{\mathrm{Ch}(R\text{-Mod})}^{\quad1}(_RZ^{\bullet},~_RA^{\bullet})=0$
as $\Ext_R^{1}(_RZ^{k},~_RA^{l})=0$ for all $k,l\in \mathbb{Z}$, and hence $_RZ^{\bullet}\in ~_R\mathcal{W}^{\bullet}$.
By similar arguments that we have used, one can see that, in the short exact sequence $(\dag)$,
if any two of $_RX^{\bullet}$, $_RY^{\bullet}$ and $_RZ^{\bullet}$  belong to $_R\mathcal{W}^{\bullet}$ then so does the third.
This shows that $_R\mathcal{W}^{\bullet}$ is thick.
Consequently, $(_R\mathcal{W}^{\bullet},~\Ch(_R\mathcal{A}\text{-}\Inj))$ forms an injective cotorsion pair
by~\cite[Proposition 3.6]{Gil2016} and by noting that $\Ch(_R\mathcal{A}\text{-}\Inj)$ contains all injective complexes.

Now we prove that
there is a cotorsion triple $$(\dg(^{\perp}{_R\mathcal{A}})\cap \Ch(_R\mathcal{A}),~\widetilde{_R\mathcal{A}},~\Ch(_R\mathcal{A}\text{-}\Inj))$$
in the exact category $\Ch(_R\mathcal{A})$ such that the left-hand cotorsion pair is projective  and the right-hand one is injective.

Indeed, the left-hand pair $(\dg(^{\perp}{_R\mathcal{A}})\cap \Ch(_R\mathcal{A}), ~\widetilde{_R\mathcal{A}})$
forms a  projective cotorsion pair due to \cite[dual of Proposition 7.3]{Gil16}. To obtain the required injective cotorsion pair $(\widetilde{_R\mathcal{A}},~\Ch(_R\mathcal{A}\text{-}\Inj))$,
by \cite[dual of Proposition 7.2]{Gil16}, it suffices to show that $$_R\mathcal{W}^{\bullet}\cap \Ch(_R\mathcal{A})=~\widetilde{_R\mathcal{A}},$$
where $_R\mathcal{W}^{\bullet}$ is as above. Using Fact 1, one can check that this holds, by using the same argument as in the proof of \cite[Theorem 4.8(4)]{Gil201703}.
\end{proof}

\begin{lem} \label{7-0}Let $R$ be a flat-typed $(\mathcal{L},\mathcal{A})$-Gorenstein ring with respect to a bi-complete
duality pair $(\mathcal{L},\mathcal{A})$. Then
\begin{enumerate}
\item There is a triangle-equivalence
$\K(R\text{-}\GProj)\simeq \K(R\text{-}\GInj)$
which restricts to a triangle equivalence $\K(R\text{-}\Proj)\simeq \K(R\text{-}\Inj)$.
\item $_R\mathcal{F}^\wedge$ (resp., ${\mathcal{F}_R}^\wedge$) is definable.
\end{enumerate}
\end{lem}
\begin{proof}(1) For any bi-complete duality pair $(\mathcal{L}, \mathcal{A})$ and any flat-typed $(\mathcal{L}, \mathcal{A})$-Gorenstein ring $R$,
by definition, all modules in $_R\mathcal{A}$ and $\mathcal{A}_R$ have finite flat dimension at most $m$ for some nonnegative integer $m$.
It follows from \cite[Proposition 2.3]{Gil2019} that all injective left and right $R$-modules have finite flat dimension at most $m$ for some nonnegative integer $m$.
Now \cite[Theorem 5.3]{Emm2012} shows that $R$ admits $\Gwgldim(R)\leq m$, and then Theorem \ref{8} applies.

(2) We only prove the case for $_R\mathcal{F}^\wedge$ since the case for ${\mathcal{F}_R}^\wedge$ can be done in a similar way.
Note from (1) that any flat-typed $(\mathcal{L}, \mathcal{A})$-Gorenstein ring $R$ admits $\Gwgldim(R)\leq m$ for some nonnegative integer $m$.
Now, by \cite[Lemma 2.1]{WZ2023} and Theorem \ref{3-3}, there is a complete and hereditary cotorsion pair
$(_R\mathcal{GP},~_R\mathcal{F}^\wedge)$ in $R\text{-}\Mod$.
This shows that $_R\mathcal{F}^\wedge$ is closed under direct products.
In addition, it is well known that $_R\mathcal{F}^\wedge$ is closed under direct limits.
Finally, $_R\mathcal{F}^\wedge$ is closed under pure submodules by \cite[Lemma 9.1.4]{EJ2000}.
\end{proof}

The next result considers some triangle-equivalences involving  $\K(_R\mathcal{A}\text{-}\Inj)$.

\begin{prop} \label{7-1}Let $R$ be a flat-typed $(\mathcal{L},\mathcal{A})$-Gorenstein ring with respect to a bi-complete
duality pair $(\mathcal{L},\mathcal{A})$, and let $\K(_R\mathcal{A}\text{-}\Inj)$ be as above.
Then we have triangle-equivalences
$$\K(R\text{-}\Inj)=\K(_R\mathcal{A}\text{-}\Inj)\simeq \K(^{\perp}{_R\mathcal{A}}\cap{_R\mathcal{A}}).$$
 \end{prop}
\begin{proof}

Let us first show the equality. Indeed, it is obvious that $\K(_R\mathcal{A}\text{-}\Inj)\subseteq \K(R\text{-}\Inj)$.
Conversely, let $_RX^{\bullet}\in \K(R\text{-}\Inj)$.
For any complex $_RA^{\bullet}\in \widetilde{_R\mathcal{A}}$, one has $_RA^{\bullet}$ is an exact complex with each cycle module $Z^{k}(_RA^{\bullet})$ in $_R\mathcal{A}$.
Note from \cite[Proposition 4.5]{WD2020} that any module in $_R\mathcal{A}$ has finite flat dimension.
Hence $_RA^{\bullet}$ has finite flat dimension in \Ch$(R\text{-Mod})$.
 So applying \cite[Corollary 3.3]{Gil201702} to the injective cotorsion pair $(^{\perp}\Ch(R\text{-Inj}),\Ch(R\text{-Inj)})$ in \Ch$(R\text{-Mod})$, one gets that $_RA^{\bullet}\in ~^{\perp}\Ch(R\text{-Inj})$.
Thus, any morphism $f^{\bullet}:~_RA^{\bullet}\to ~_RX^{\bullet}$ with $_RA^{\bullet}\in \widetilde{_R\mathcal{A}}$ is null-homotopic.
This shows that $_RX^{\bullet}\in\K(_R\mathcal{A}\text{-}\Inj)$.

Now we prove the remaining triangle-equivalence. By Lemma \ref{7}, it suffices to see an equality $\Ch(^{\perp}{_R\mathcal{A}})=\dg(^{\perp}{_R\mathcal{A}})$.
Note that the subcategory $_R\mathcal{A}$ is closed under pure submodules and pure quotients of modules.
Hence, there is a set $\{_RA_{\lambda}~|~\lambda\in \Lambda\}\subseteq ~_R\mathcal{A}$ such that
any left $R$-module in $_R\mathcal{A}$ is a transfinite extension of   $\{_RA_{\lambda}~|~\lambda\in \Lambda\}$.
We put $L_R=\coprod_{\lambda\in \Lambda}(_RA_{\lambda})^{+}$.
Then $(L_R)^{+}\cong\prod_{\lambda\in \Lambda}(_RA_{\lambda})^{++}$.
Note further that $\coprod_{\lambda\in \Lambda} {_R} A_{\lambda}$ (resp. each $_RA_{\lambda}$)
is a pure submodule of $\coprod_{\lambda\in \Lambda}(_RA_{\lambda})^{++}$ (resp. $(_RA_{\lambda})^{++}$)
and that any transfinite extension is a special direct limit.
It is then easy to see that $_R\mathcal{A}=\overline{\{L_R^{+}\}}$.
Now let $I_R=L_R$ and $$0\to~ _RM\to~_RP\to ~_RM \to0$$ be a short exact sequence of left $R$-modules with $_RP\in~^{^{\perp}}\overline{\{I_R^{+}\}}$.
Due to \cite[Proposition 4.5]{WD2020}, $I_R$ has finite flat dimension, and hence, $(I_R)^{J}$ has finite flat dimension for any set $J$ by Lemma \ref{7-0}(2).
Since $\Ext_R^{i\geq1}(((I_R)^{J})^{+}, ~_RP)=0$, using the isomorphism
$$\Tor^R_{i\geq1}((I_R)^{J},~_RP)^{+}\cong \Ext_R^{i\geq1}(((I_R)^{J})^{+}, ~_RP)$$
one has $\Tor ^R_{i\geq1}((I_R)^{J},~_RP)=0$.
It is now a routine to check that $\Tor^R_{i\geq1}((I_R)^{J},~_RM)=0$.
Whence, one can use \cite[Proposition 4.2]{SS2020} to obtain that $_RM\in ~^{\perp}{\overline{\{L_R^{+}\}}}$, i.e., $_RM\in~^{\perp}{_R\mathcal{\mathcal{A}}}$.
Thus, by a dual of \cite[Theorem 5.3 and its proof]{BCE2020}, one can conclude that $\Ch(^{\perp}{_R\mathcal{\mathcal{A}}})=\dg(^{\perp}{_R\mathcal{\mathcal{A}}})$, as required.
\end{proof}

Now combining with Lemma \ref{7-0}(1), Propositions \ref{7-4} and \ref{7-1}, we have

\begin{thm} \label{7-5}Let $R$ be a flat-typed $(\mathcal{L},\mathcal{A})$-Gorenstein ring with respect to a bi-complete
duality pair $(\mathcal{L},\mathcal{A})$.
Then we have a triangle-equivalence
$$\K({_R\mathcal{F}}\cap~_R\mathcal{C}ot) \simeq \K(^{\perp}{_R\mathcal{A}}\cap{_R\mathcal{A}}).$$
 \end{thm}

As mentioned above, any Gorenstein $n$-coherent ring $R$ admits $\Gwgldim(R)<\infty$.
On the other hand, it is known that the pair $(\mathcal{FP}_n\mathcal{F}, \mathcal{FP}_n\mathcal{I})$ forms a bi-complete duality if and only if the ground ring is two-sided $n$-coherent (see \cite[Examples 3.7-3.8]{WD2020}).
Furthermore, a Gorenstein $n$-coherent ring is just a flat-typed $(\mathcal{FP}_n\mathcal{F}, \mathcal{FP}_n\mathcal{I})$-Gorenstein ring (see \cite[Example 4.3]{WD2020}).
In particular, Ding-Chen rings (i.e., Gorenstein $1$-coherent rings) are just flat-typed $(\mathcal{F}, \mathcal{FI})$-Gorenstein rings.
Thus, as an immediate consequence of Theorem \ref{7-5}, we have

\begin{cor} \label{7-5-1} Let $n\in \mathbb{N}\cup \{\infty\}$ and $R$ be a Gorenstein $n$-coherent ring.
Then there is a triangle-equivalence
$$\K({_R\mathcal{F}}\cap~_R\mathcal{C}ot)\simeq \K(^{\perp}{_R\mathcal{FP}_n\mathcal{I}}\cap{_R\mathcal{FP}_n\mathcal{I}}).$$
Particularly, if $R$ is a Ding-Chen ring, then there exists a triangle-equivalence
$$\K({_R\mathcal{F}}\cap~_R\mathcal{C}ot)\simeq \K({_R\mathcal{FI}}\cap~_R\mathcal{FP})$$
in which $_R\mathcal{FP}=~{^{\perp}{_R\mathcal{FI}}}$ means the subcategory of $R$-Mod consisting of FP-projective left $R$-modules.
 \end{cor}

Given a left $R$-module $_RM$ and an integer $k\in \mathbb{Z}$, following \cite{Gil201702}, we denote by $\mathrm{S}^k(_RM)$ the complex $\cdots\to 0 \to~_RM \to 0 \to\cdots$ with
$_RM$ in the $k$th degree.  Usually, $\mathrm{S}^0(_RM)$ is denoted
simply by $_RM$.

Let $\mathcal{T}$ be a triangulated category with small coproducts. Recall that an object $C$ of $\mathcal{G}$ is \emph{compact} if
for each collection $\{T_j~|~j\in J\}$ of objects of $\mathcal{T}$, the canonical morphism
$$\coprod_{j\in J}\Hom_{\mathcal{T}}(C,T_j)\to \Hom_{\mathcal{T}}(C,\coprod_{j\in J}T_j)$$
is an isomorphism. The category $\mathcal{T}$ is \emph{compactly generated} if there exists a small set $\mathcal{U}\subseteq \mathcal{G}$ of compact
objects such that for each $0\neq T\in\mathcal{T}$ there is a morphism $0 \neq f : \Sigma^{m}U \to T$ for some $U\in \mathcal{U}$ and $m\in \mathbb{Z}$, where $\Sigma$ denotes the autofunctor of $\mathcal{T}$.

\begin{lem} \label{7-6}Let $R$ be a ring with a bi-complete
duality pair $(\mathcal{L},\mathcal{A})$.
If the cogenerating set $_R\mathcal{S}$ of the cotorsion pair $(^{\perp}{_R\mathcal{A}},~_R\mathcal{A})$ consists of left $R$-modules of type FP$_\infty$ and contains the regular module $_RR$,
then $\K(_R\mathcal{A}\text{-}\Inj)$ is compactly generated.
 \end{lem}
\begin{proof} We consider the set $$_R\mathcal{V}^{\bullet}=\{_RV^{\bullet}=\mathrm{S}^{k}(_RM)~|~_RM\in~_R\mathcal{S}, k\in \mathbb{Z}\}$$ induced by $_R\mathcal{S}$.
Then \cite[Lemma 2.3]{Gil2017} yields that every complex in $_R\mathcal{V}$ is compact.

To see that $\K(_R\mathcal{A}\text{-}\Inj)$ is generated by $_R\mathcal{V}^{\bullet}$, we let $_RX^{\bullet}$ be a complex in $\K(_R\mathcal{A}\text{-}\Inj)$ such that
$$\Hom_{\mathrm{K}(_R\mathcal{A}\text{-Inj})}(\Sigma^{k}(_RV^{\bullet}), ~_RX^{\bullet})=0$$ for all $k\in \mathbb{Z}$ and all $_RV^{\bullet}\in ~_R\mathcal{V}^{\bullet}$.
Since $_R\mathcal{S}$ contains $_RR$, one has that $_RX^{\bullet}$ is exact. Now, by \cite[Theorem 3.5]{YgL2011}
there exists a short exact sequence in \Ch$(R\text{-Mod})$
$$(\ddag)\quad0\to ~_RX^{\bullet}\to  ~_RA^{\bullet}\to ~_RT^{\bullet}\to0$$
with $_RA^{\bullet}\in ~\dg(_R\mathcal{A})$ and $_RT^{\bullet}\in ~\widetilde{^{\perp}{_R\mathcal{A}}}$.
Note that $_RA^{\bullet}$ is exact since $_RX^{\bullet}$ and $_RT^{\bullet}$ are so.
Hence, one gets that $_RA^{\bullet}\in ~\widetilde{_R\mathcal{A}}$ by \cite[Theorem 3.12]{Gil2004},
which implies that $$\Hom_{\mathrm{K}(R\text{-}\mathrm{Mod})}(\Sigma^{k}(_RV^{\bullet}), ~_RA^{\bullet})=0$$
for all $k\in \mathbb{Z}$ and all $_RV^{\bullet}\in ~_R\mathcal{V}^{\bullet}$.
Now one can obtain that
 $$\Hom_{\mathrm{K}(R\text{-}\mathrm{Mod})}(\Sigma^{k} (_RV^{\bullet}), ~_RT^{\bullet})=0$$
for all $k\in \mathbb{Z}$ and all $_RV^{\bullet}\in ~_R\mathcal{V}^{\bullet}$.
On the other hand, one has that each component module of $_RT^{\bullet}$ belongs to $_R\mathcal{A}$ since so does $_RX^{\bullet}$ and $_RA^{\bullet}$. It follows that
each component module of $_RT^{\bullet}$ belongs to $^{\perp}{_R\mathcal{A}}\cap{_R\mathcal{A}}$, and so $\Ext_R^{1}(_RM, ~_RT^{k})=0$
for all $k\in \mathbb{Z}$ and all $~_RM\in~_R\mathcal{S}\subseteq ~^{\perp}{_R\mathcal{A}}$.
Whence,
 $$\Ext_R^{1}(_RM, ~\Z^{k}(_RT^{\bullet}))\cong \Ext_{\mathrm{Ch}(R\text{-Mod})}^{\quad1}(\mathrm{S}^{k}(_RM), ~_RT^{\bullet})=\Hom_{\mathrm{K}(R\text{-}\mathrm{Mod})}(\mathrm{S}^{k}(_RM), ~_RT^{\bullet})=0$$
for all $_RM\in~_R\mathcal{S}$ and all $k\in \mathbb{Z}$. It follows that $\Z^{k}(_RT^{\bullet})\in~ _R\mathcal{A}$. Now $\Z^{k}(_RT^{\bullet})\in~^{\perp}{_R\mathcal{A}}\cap{_R\mathcal{A}}$ since $_RT^{\bullet}\in ~\widetilde{^{\perp}{_R\mathcal{A}}}$.
Thus, $_RT^{\bullet}$ is a contractible complex, which is a zero object in $\K(\widetilde{^{\perp}{_R\mathcal{A}}}\cap \widetilde{_R\mathcal{A}})$, of course is a zero object in $\K(\dg(^{\perp}{_R\mathcal{A}})\cap \Ch(_R\mathcal{A}))$, and hence is a zero object in $\K(_R\mathcal{A}\text{-}\Inj)$ by Lemma \ref{7}.
Notice further that the above short exact sequence $(\ddag)$ is degree-wise split as each component module of $_RX^{\bullet}$ is injective. So, we have an isomorphism $_RX^{\bullet}\cong  ~_RA^{\bullet}$ in $\K(_R\mathcal{A}\text{-}\Inj)$.
Consequently, any morphism $f^{\bullet}:~_RX^{\bullet}\to ~_RX^{\bullet}$ is null-homotopic since $_RX^{\bullet}\cong ~_RA^{\bullet}\in ~\widetilde{_R\mathcal{A}}$. It follows that $_RX^{\bullet}$ is a contractible complex, which is a zero object in $\mathrm{K}(_R\mathcal{A}\text{-}\Inj)$.
Therefore, $\K(_R\mathcal{A}\text{-}\Inj)$ is generated by the set $_R\mathcal{V}^{\bullet}$ of compact objects.
\end{proof}

\begin{rem} \label{7-6-0} {\rm For the special duality pair $(\mathcal{FP}_\infty\mathcal{F}, \mathcal{FP}_\infty\mathcal{I})$, the corresponding result in Lemma \ref{7-6} was proved in \cite[Theorem 4.4]{Gil2017}
in terms of the language of model structures. Compared with it,  the proof of Lemma \ref{7-6} is more direct.}
\end{rem}

\begin{thm} \label{7-7}Let $R$ be a flat-typed $(\mathcal{L},\mathcal{A})$-Gorenstein ring with respect to a bi-complete
duality pair $(\mathcal{L},\mathcal{A})$. If the cogenerating set $_R\mathcal{S}$ of the cotorsion pair $(^{\perp}{_R\mathcal{A}},~_R\mathcal{A})$ consists of left $R$-modules of type FP$_\infty$ and contains the regular module $_RR$,
then
$$\K({_R\mathcal{F}}\cap~_R\mathcal{C}ot)\simeq \K(_R\mathcal{L}\text{-}\Proj)=\K(R\text{-}\Proj)\simeq\K(R\text{-}\Inj)=\K(_R\mathcal{A}\text{-}\Inj)\simeq \K(^{\perp}{_R\mathcal{A}}\cap{_R\mathcal{A}})$$
are compactly generated.
 \end{thm}
\begin{proof} Note that all the categories are triangle equivalent by Lemma \ref{7-0}(1), Propositions \ref{7-1} and \ref{7-4}. Thus, all of them are compactly generated by Lemma \ref{7-6}.
\end{proof}

We end this section with an application of Theorem \ref{7-7}.

\begin{cor} \label{7-7-0} Let $n\in \mathbb{N}\cup \{\infty\}$ and $R$ be a Gorenstein $n$-coherent ring.
Then
$$\K({_R\mathcal{F}}\cap~_R\mathcal{C}ot))\simeq \K(R\text{-}\Proj)\simeq\K(R\text{-}\Inj)\simeq \K(^{\perp}{_R\mathcal{FP}_n\mathcal{I}}\cap{_R\mathcal{FP}_n\mathcal{I}})$$
are compactly generated.
 \end{cor}
\begin{proof} Note from \cite[Example 4.3]{WD2020} that $R$ is just a flat-typed $(\mathcal{FP}_n\mathcal{F}, \mathcal{FP}_n\mathcal{I})$-Gorenstein ring.
In addition, the cogenerating set $_R\mathcal{S}$ of the cotorsion pair $(^{\perp}{_R\mathcal{FP}_n\mathcal{I}},~_R\mathcal{FP}_n\mathcal{I})$  consists of
finitely $n$-presented left $R$-modules which coincide with left $R$-modules of type FP$_\infty$ as $R$ is two-sided $n$-coherent (see \cite[Proposition 4.1 and Corollary 4.2 as well as Theorem 2.4]{BP2017}).
Thus, Theorem \ref{7-7} applies.
\end{proof}

\begin{rem} \label{7-7-2} {\rm Let $n\in \mathbb{N}\cup \{\infty\}$ and $R$ be a Gorenstein $n$-coherent ring. The fact that $\K(R\text{-}\Inj)$ is compactly generated can be also obtained by \cite[Theorem 5.3]{WLY2018} and  \cite[Theorem 4.4]{Gil2017}.}
\end{rem}

\section{\bf Grothendieck duality}

Let $S$ be a
left coherent ring and $R$ be a right coherent ring.
A complex $_RD_S^{^{\bullet}}$ of $R$-$S$-bimodules
is said to be a \emph{dualizing complex} \cite{Pos2017} for the rings $R$ and $S$ if the following
conditions are satisfied:
\begin{enumerate}
\item[(DC1)] all component modules of $_RD_S^{^{\bullet}}$ are FP-injective left $R$-modules and FP-injective right
$S$-modules;
\item[(DC2)] all homology modules of $_RD_S^{^{\bullet}}$ are finitely presented left $R$-modules and finitely presented right $S$-modules;
\item[(DC3)] the homothety maps $R\to \mathcal{H}om_{S^{op}}(_RD_S^{\bullet}, ~_RD_S^{\bullet})$
) and $S\to \mathcal{H}om_{R}(_RD_S^{\bullet}, ~_RD_S^{\bullet})$
) are
quasi-isomorphisms of DG-rings.
\end{enumerate}
This notion coincides with the notion of ``dualizing complex'' in \cite[Definition 1.1]{CFH2006} when $R$ and $S$ are left and right noetherian, respectively.
In addition, we say that a left coherent ring $R$ \emph{has a dualizing complex} $_RD_R^{\bullet}$ if $_RD_R^{\bullet}$ is a dualizing complex for the rings $R$ and $R^{op}$.

Grothendieck duality is a classical subject which can go back 1958s. Roughly speaking, it is a statement concerning the existence of a right adjoint to
the ``direct image with compact support'' functor between derived categories of sheaves or modules.
For a ring $R$, over the years many people investigated Grothendieck duality for derived categories of $R$-modules by providing the following insights:
\begin{enumerate}
\item[(\textbf{GD1})] $R$ has a duality complex.
\item[(\textbf{GD2})] There is a triangle-equivalence $\mathbf{D}^{b}(R^{op}\text{-}\mathrm{mod})^{op}\simeq \mathbf{D}^{b}(R\text{-}\mathrm{mod}),$
where the right category is the bounded derived category of finitely presented $R$-modules and the left one is the opposite category of the bounded derived category of finitely  presented right $R$-modules;
  \item[(\textbf{GD3})] There is a triangulated equivalence $\K(R\text{-}\p)\simeq \K(R\text{-}\inj)$.
 \end{enumerate}
According to \cite[Introduction]{Nee2008}, we know that, for a left noetherian and right coherent ring such that all flat left $R$-modules have finite projective dimension, there are implications
$$(\textbf{GD1})\Longrightarrow (\textbf{GD3}) \Longrightarrow(\textbf{GD2}).$$

As shown in the next remark, the condition on the ring such that the above implications
hold can be relaxed.

\begin{rem} \label{10-12} {\rm Let $R$ be a two-sided coherent ring
and let $(\textbf{GD1})$, $(\textbf{GD2})$ and $(\textbf{GD3})$ be as above.
Then $(\textbf{GD3}) \Longrightarrow(\textbf{GD2})$ follows from \cite[Corollary 9(2)]{WLY2021}.
On the other hand, the implication $(\textbf{GD1})\Longrightarrow (\textbf{GD3})$ holds if all FP-injective left $R$-modules have finite
injective dimension. Indeed, suppose that $R$ has a dualizing complex. Using \cite[Proposition 4.3]{Pos2017} one gets that all flat right $R$-modules have finite projective dimension.
So we obtain that $\K(R\text{-}\p)\simeq \K(R\text{-}\inj)$ by applying \cite[Theorems 2.4, 3.2 and 4.4]{Pos2017} as $R$ is a a two-sided coherent ring such that all FP-injective left $R$-modules have finite
injective dimension. However, in order that $(\textbf{GD1})\Longrightarrow (\textbf{GD3})$ holds, we do not know whether the condition that all FP-injective left $R$-modules have finite
injective dimension can be omitted.}
\end{rem}

Recall that the notation $_R\mathcal{FP}_\infty\mathcal{F}$ (resp $_R\mathcal{FP}_\infty\mathcal{I}$) stands for the subcategory of $R$-Mod
consisting of all level (resp. absolutely clean) modules. For simplicity, we denote them by $_R\mathcal{L}e$ and $_R\mathcal{A}c$, respectively.
Note that, over any ring $R$, the pair $(\mathcal{L}e, \mathcal{A}c)$ is a bi-complete duality pair (see \cite[Example 3.6]{WD2020}).
Recall from \cite[Definition 5.1]{Gil2017} that a complex $_RA^{\bullet}$ is \emph{AC-injective} if $_RA^{\bullet}$ is a complex in $\Ch(_R\mathcal{A}c\text{-}\Inj)$.
The corresponding homotopy subcategory consisting of all AC-injective complexes will be denoted by $\K(_R\mathcal{A}c\text{-}\Inj)$.

Following \cite{WLY2021}, we denote by $\mathbf{D}^{b}(R\text{-}\mathrm{tmod}$) (resp. $\mathbf{D}^{b}(R^{op}\text{-}\mathrm{tmod})^{op}$)
the bounded derived category of left $R$-modules of type FP$_\infty$ (resp. the opposite category of the bounded derived category of right $R$-modules of type FP$_\infty$),
where $R\text{-}\mathrm{tmod}$ (resp. $R^{op}\text{-}\mathrm{tmod}$) denotes the subcategory of $R$-Mod (resp. $R^{op}$-Mod) consisting of all modules of type FP$_\infty$.

Now we can give another application of Theorem \ref{8} (or Corollary \ref{9-2}), which can be compared with Remark \ref{10-12}.

\begin{cor} \label{11}
Let $R$ be a ring such that $\K(R\text{-}\p)$ is compactly generated and that there is an equality $\K(_R\mathcal{A}c\text{-}\Inj)=\K(R\text{-}\inj)$. Consider the following conditions
\begin{enumerate}
\item $R$ is a Gorenstein $n$-coherent ring, where $n\in \mathbb{N}\cup\{\infty\}$;
\item There is a triangle-equivalence $\mathbf{D}^{b}(R^{op}\text{-}\mathrm{tmod})^{op}\simeq \mathbf{D}^{b}(R\text{-}\mathrm{tmod});$
\item There is a triangle-equivalence $\K(R\text{-}\p)\simeq \K(R\text{-}\inj)$.
 \end{enumerate}
Then we have the implications $(1)\Rightarrow(3)\Rightarrow(2)$.
\end{cor}

\begin{proof} $(1)\Rightarrow(3)$ follows from Corollary \ref{9-2}.
 $(3)\Rightarrow(2)$ comes from \cite[Corollary 9(1)]{WLY2021}.
\end{proof}

Recall that $\mathbf{D}^{b}(R\text{-}\mathrm{mod}$) (resp. $\mathbf{D}^{b}(R^{op}\text{-}\mathrm{mod})^{op}$) stands for
the bounded derived category of finitely presented left $R$-modules (resp. the opposite category of the bounded derived category of finitely presented right $R$-modules).
Let us consider the particular case of $n=1$ in Corollary \ref{11}.

\begin{cor} \label{10}
Let $R$ be a left and right coherent ring. Consider the following conditions

\begin{enumerate}
\item $R$ is a Ding-Chen ring;
\item There is a triangle-equivalence $\mathbf{D}^{b}(R^{op}\text{-}\mathrm{mod})^{op}\simeq \mathbf{D}^{b}(R\text{-}\mathrm{mod});$
\item There is a triangle-equivalence $\K(R\text{-}\p)\simeq \K(R\text{-}\inj)$.
 \end{enumerate}
Then we have the implications $(1)\Rightarrow(3)\Rightarrow(2)$.
 \end{cor}

\begin{proof}
Suppose that $R$ is left and right coherent. Then $\K(R\text{-}\p)$ is compactly generated and that there is an equality $\K(_R\mathcal{A}c\text{-}\Inj)=\K(R\text{-}\inj)$ (see \cite[Proposition 7.14]{Nee2008} and \cite[Remark in p.109]{Gil2017}). In addtion, we know from \cite[Theorem 3.21]{Gil201703} that the two subcategories $R^{op}\text{-mod}$ and $R^{op}\text{-tmod}$ (resp. $R\text{-mod}$ and $R\text{-tmod}$)
coincide. Now Corollary \ref{11} applies as Ding-Chen rings are exactly Gorenstein $1$-coherent rings.
\end{proof}

\begin{rem} \label{10-10} {\rm We have a ``left Gorenstein version'' of Corollary \ref{10}, that is, there are implications $(1)\Rightarrow(3)\Rightarrow(2)$
among the following conditions for any left and right coherent ring $R$:
\begin{enumerate}
\item $R$ is a left Gorenstein ring, that is, $R$ is a ring with $\Ggldim(R)<\infty$;
\item There is a triangle-equivalence $\mathbf{D}^{b}(R^{op}\text{-}\mathrm{mod})^{op}\simeq \mathbf{D}^{b}(R\text{-}\mathrm{mod});$
\item $\K(R\text{-}\p)$ and $\K(R\text{-}\inj)$ are triangle equivalent.
 \end{enumerate}
Indeed,  $(1)\Rightarrow(3)$ and $(3)\Rightarrow(2)$ follow from \cite[Theorem B]{Chen2010} and \cite[Corollary 9(1)]{WLY2021}, respectively.}
\end{rem}

For any triangulated category $\mathcal{T}$, denote by $\mathcal{T}^{c}$ the triangulated subcategory of $\mathcal{T}$ consisting of all compact objects.
It is well known that the usual (unbounded) derived category $\mathbf{D}(R\text{-Mod})$ is compactly generated with a triangle equivalence
$\mathbf{D}(R\text{-}\Mod)^{c}\simeq \K^b(R\text{-proj}),$ where $\K^{b}(R\text{-}\mathrm{proj})$ means the bounded homotopy subcategory of
 finitely presented (equivalently finitely generated) projective left $R$-modules. The next lemma can be compared with such a fact, where
$\K^{b}(R^{op}\text{-}\mathrm{proj})^{op}$ denotes the opposite category of the bounded homotopy subcategory of
 finitely presented (equivalently finitely generated) projective right $R$-modules.

\begin{lem} \label{10-11} Let $R$ be a ring. Then there is a triangle-equivalence
$$\mathbf{D}(R\text{-}\Mod)^{c}\simeq \K^{b}(R^{op}\text{-}\mathrm{proj})^{op}.$$
\end{lem}
\begin{proof}
For any  $X^{\bullet}\in \K(R\text{-}\mathrm{proj})$, we define
$${X^{\bullet}}^{*}=\Hom_R(X^{\bullet},R)\quad\text{and}\quad {X^{\bullet}}^{**}=\Hom_{R^{op}}(\Hom_R(X^{\bullet},R),R).$$

It is a routine to check that the canonical morphisms of complexes $X^{\bullet}\to {X^{\bullet}}^{**}$ and ${X^{\bullet}}^{*}\to {X^{\bullet}}^{***}$ are isomorphisms for all $X^{\bullet}\in \K(R\text{-}\mathrm{proj})$.
It follows that
$$\xymatrix@C=95pt{
\K(R\text{-}\mathrm{proj}) \ar@<0.4ex>[r]^{(-)^{*}} & \K(R^{op}\text{-}\mathrm{proj})^{op}
   \ar@<0.4ex>[l]^{(-)^{*}}.}$$
are quasi-inverse equivalences of triangulated categories.

It is well-known that $\mathbf{D}(R\text{-}\Mod)^{c}\simeq \K^b(R\text{-}\mathrm{proj})$, now the result follows.
\end{proof}

Let $\mathcal{T}, \mathcal{T}', \mathcal{T}''$ be triangulated categories and $F: \mathcal{T}'\to \mathcal{T}$, $G: \mathcal{T}\to \mathcal{T}''$ triangulated functors.
Recall that $\mathcal{T}'\overset{F}\to \mathcal{T}\overset{G}\to \mathcal{T}''$ is
a \emph{localization sequence} if the following holds:\begin{enumerate}
  \item[(LS1)] The functor $F$ has a right adjoint $F_\rho : \mathcal{T} \to \mathcal{T}'$ satisfying $F_\rho F = \mathrm{Id}_{{\mathcal{T}}'} $.
\item[(LS2)] The functor $G$ has a right adjoint $G_\rho : \mathcal{T}'' \to \mathcal{T}$ satisfying $G_\rho G = \mathrm{Id}_{{\mathcal{T}}''} $.
\item[(LS3)] Let $T$ be an object in $\mathcal{T}$. Then $G(T) = 0$ if and only if $T = F(T')$ for some $T' \in \mathcal{T}'$.
\end{enumerate}

Let $\mathcal{T}'\overset{F}\to \mathcal{T}\overset{G}\to \mathcal{T}''$ be a localization sequence of triangulated categories. It is said to be a \emph{recollement} if it is also a colocalization sequence.
Here we say that $\mathcal{T}'\overset{F}\to \mathcal{T}\overset{G}\to \mathcal{T}''$  is a \emph{colocalization sequence} if it satisfies the following conditions:\begin{enumerate}
  \item[(CLS1)] The functor $F$ has a left adjoint $F_\rho : \mathcal{T} \to \mathrm{T}'$ satisfying $F_\rho F = \mathrm{Id}_{{\mathcal{T}}'} $.
\item[(CLS2)] The functor $G$ has a left adjoint $G_\rho : \mathcal{T}'' \to \mathrm{T}$ satisfying $G_\rho G = \mathrm{Id}_{{\mathcal{T}}''} $.
\item[(CLS3)] Let $T$ be an object in $\mathcal{T}$. Then $G(T) = 0$ if and only if $T = F(T')$ for some $T' \in \mathcal{T}'$.
\end{enumerate}

\begin{rem} \label{10-13} {\rm As remarked in \cite{BI2007}, the notions of torsion pairs and recollements are closely related.}
\end{rem}

Note that any finitely generated (finitely presented ) projective module is of type FP$_\infty$.
This implies that $\K^{b}(R\text{-proj})$ (resp. $\K^{b}(R^{op}\text{-proj})$)
is a thick subcategory of $\mathbf{D}^{b}(R\text{-tmod})$ (resp. $\mathbf{D}^{b}(R^{op}\text{-tmod})$).
And hence we have the corresponding Verdier quotients.

\begin{thm} \label{12}
Let $R$ be a ring such that $\K(R\text{-}\p)$ is compactly generated and that there is an equality $\K(_R\mathcal{A}c\text{-}\Inj)=\K(R\text{-}\inj)$. Consider the following conditions
 \begin{enumerate}
\item $R$ is a Gorenstein $n$-coherent ring, where $n\in \mathbb{N}\cup\{\infty\}$ or $R$ is left Gorenstein (that is, G$\text{-}$gldim$(R)$ $<\infty$);
\item  There is a triangle-equivalence $\mathbf{D}^{b}(R^{op}\text{-}\mathrm{tmod})/\K^{b}(R^{op}\text{-}\mathrm{proj})^{op}\simeq \mathbf{D}^{b}(R\text{-}\mathrm{tmod})/\K^{b}(R\text{-}\mathrm{proj}),$
where the left category is the opposite category of the corresponding Verdier quotient triangulated category of right $R$-modules;
  \item $\K_{\mathrm{ac}}(R\text{-}\p)$ and $\K_{\mathrm{ac}}(R\text{-}\inj)$ are triangle equivalent.
 \end{enumerate}
Then one has implications
$(1)\Longrightarrow (3) \Longrightarrow(2).$
\end{thm}
\begin{proof}
$(1)\Rightarrow(3)$ follows from Corollaries \ref{9-0-0}(1) and \ref{9-2}.

$(3)\Rightarrow(2)$. Suppose that there is a triangle equivalence $\K_{\mathrm{ac}}(R\text{-}\Proj)\simeq\K_{\mathrm{ac}}(R\text{-}\Inj)$.
According to \cite[Corollaries 2.2.2 and 2.2.3]{Bec2014}, we have two recollements
\[
\xymatrixrowsep{5pc}
\xymatrixcolsep{6pc}
\xymatrix
{\K_{\ac}(R\text{-}\inj) \ar[r]^{\mathrm{I(inj)}} &\K(R\text{-}\inj) \ar@<-3ex>[l]_{\mathrm{I(inj)}^{l}} \ar@<3ex>[l]^{\mathrm{I(inj)}_{r}} \ar[r]^{\mathrm{Q(inj)}}  &\mathbf{D}(R\text{-}\Mod) \ar@<-3ex>[l]_{\mathrm{Q(inj)}^{l}} \ar@<3ex>[l]^{\mathrm{Q(inj)}_{r}}
  } \]
and
\[
\xymatrixrowsep{5pc}
\xymatrixcolsep{6pc}
\xymatrix
{\K_{\ac}(R\text{-}\p) \ar[r]^{\mathrm{I(prj)}} &\K(R\text{-}\p) \ar@<-3ex>[l]_{\mathrm{I(prj)}^{l}} \ar@<3ex>[l]^{\mathrm{I(prj)}_{r}} \ar[r]^{\mathrm{Q(proj)}}  &\mathbf{D}(R\text{-}\Mod) \ar@<-3ex>[l]_{\mathrm{Q(proj)}^{l}} \ar@<3ex>[l]^{\mathrm{Q(proj)}_{r}}
  .} \]
In particular, we have two colocalization sequences
$$\xymatrix@C=0.5cm{
 & \mathbf{D}(R\text{-}\Mod) \ar[r]^{\mathrm{I(inj)}^{l}} & \K(R\text{-}\inj) \ar[r]^{\mathrm{Q(inj)}^{l}} &\K_{\ac}(R\text{-}\inj) }\quad\text{and}~\xymatrix@C=0.5cm{
 & \mathbf{D}(R\text{-}\Mod) \ar[r]^{\mathrm{I(prj)}^{l}} & \K(R\text{-}\p) \ar[r]^{\mathrm{Q(prj)}^{l}} &\K_{\ac}(R\text{-}\p)}.$$
It is well-known that $\mathbf{D}(R\text{-}\Mod)$ is compactly generated.
In addition, $\K(R\text{-}\inj)$ (resp. $\K(R\text{-}\p)$) is compactly generated since $R$ is left and right coherent (see \cite[Proposition 7.14]{Nee2008} and \cite[Corollary 6.13]{Sto2014}).
Now, by an argument used in the proof of \cite[Proposition 5.3]{Kra2005}, one has $\K_{\ac}(R\text{-}\inj)$ and $\K_{\ac}(R\text{-}\p)$) are compactly generated
with the following triangle-equivalences
$$\K_{\ac}(R\text{-}\inj)^{\text{c}}\simeq \K(R\text{-}\inj)^{\text{c}}/\mathbf{D}(R\text{-}\Mod)^{\text{c}}\quad\text{and}\quad \K_{\ac}(R\text{-}\p)^{\text{c}}\simeq \K(R\text{-}\p)^{\text{c}}/\mathbf{D}(R\text{-}\Mod)^{\text{c}}.$$
On one hand,  using the equality $\K(_R\mathcal{A}c\text{-}\Inj)=\K(R\text{-}\inj)$ and \cite[Theorem 8(2)]{WLY2021},
one gets that $$\K(R\text{-}\inj)^{\text{c}}\simeq \mathbf{D}^{b}(R\text{-}\mathrm{mod}).$$
Hence,
$$\K_{\ac}(R\text{-}\inj)^{\text{c}}\simeq \K(R\text{-}\inj)^{\text{c}}/\mathbf{D}(R\text{-}\Mod)^{\text{c}}=\mathbf{D}^{b}(R\text{-}\mathrm{mod})/\K^{b}(R\text{-}\mathrm{proj})$$
as $\mathbf{D}(R\text{-}\Mod)^{\text{c}}\simeq \K^{b}(R\text{-}\mathrm{proj})$.
On the other hand, using the compact generation of $\K(R\text{-}\p)$ and \cite[Theorem 8(2)]{WLY2021}, one obtains that $$\K(R\text{-}\p)^{\text{c}}\simeq\mathbf{D}^{b}(R^{op}\text{-}\mathrm{mod})^{op}.$$
Whence,
$$\K_{\ac}(R\text{-}\p)^{\text{c}}\simeq \K(R\text{-}\p)^{\text{c}}/\mathbf{D}(R\text{-}\Mod)^{\text{c}}=\mathbf{D}^{b}(R^{op}\text{-}\mathrm{mod})^{op}/\K^{b}(R^{op}\text{-}\mathrm{proj})^{op}$$
as Lemma \ref{10-11} tells us that $\mathbf{D}(R\text{-}\Mod)^{\text{c}}\simeq \K^{b}(R^{op}\text{-}\mathrm{proj})^{op}$.
Consequently, the triangle-equivalence in (2) holds by the assumption $\K_{\ac}(R\text{-}\p)\simeq \K_{\ac}(R\text{-}\inj)$ which can restrict to compact objects.
\end{proof}

\begin{rem} \label{12-0} {\rm By the proof of Theorem \ref{12}, one can see that $\K_{\ac}(R\text{-}\p)$ (resp. $\K_{\ac}(R\text{-}\inj)$) is compactly generated whenever $\K(R\text{-}\p)$ (resp. $\K(R\text{-}\inj)$) is so.}
\end{rem}

The Verdier quotient triangulated category $\mathbf{D}_{sg}(R):= \mathbf{D}^b(R\text{-}\mathrm{mod})/\K^b (R\text{-}\mathrm{proj})$ was primarily introduced and studied by
Buchweitz in \cite{Buc1986} under the name ``stable derived category'' over noetherian rings. It was renamed by Orlov in \cite{Orl2004}
and is known then as \emph{ singularity category} in the literature.
We end the paper with the next corollary, which can be viewed as a singularity category version of Corollary \ref{10} and Remark \ref{10-10}.

\begin{cor} \label{12-1}
Let $R$ be a left and right coherent ring. Consider the following conditions
\begin{enumerate}
\item[(\textbf{GD1$'$})] $R$ is Ding-Chen, or $R$ is left Gorenstein (that is $\Ggldim(R)<\infty$);
\item[(\textbf{GD2$'$})] There is a triangle-equivalence $\mathbf{D}^{b}(R^{op}\text{-}\mathrm{mod})/\K^{b}(R^{op}\text{-}\mathrm{proj})^{op}\simeq \mathbf{D}^{b}(R\text{-}\mathrm{mod})/\K^{b}(R\text{-}\mathrm{proj}),$
where the left category is the opposite category of the singularity category of right $R$-modules;
  \item[(\textbf{GD3$'$})] $\K_{\mathrm{ac}}(R\text{-}\p)$ and $\K_{\mathrm{ac}}(R\text{-}\inj)$ are triangle equivalent.
 \end{enumerate}
Then one has implications
$(\textbf{GD1$'$})\Longrightarrow (\textbf{GD3$'$}) \Longrightarrow(\textbf{GD2$'$}).$
\end{cor}
\begin{proof}
Suppose that $R$ is left and right coherent. Then $\K(R\text{-}\p)$ is compactly generated and there is an equality $\K(_R\mathcal{A}c\text{-}\Inj)=\K(R\text{-}\inj)$ (see \cite[Proposition 7.14]{Nee2008} and \cite[Remark in p.109]{Gil2017}). In addition, we know from \cite[Theorem 3.21]{Gil201703} that the two subcategories $R^{op}\text{-mod}$ and $R^{op}\text{-tmod}$ (resp. $R\text{-mod}$ and $R\text{-tmod}$)
coincide. Now the implications hold by Theorem \ref{12}.
\end{proof}

\begin{rem} \label{12-2} {\rm Using Proposition \ref{9-6-2}, one can see that $\textbf{GD3}'$ in Corollary \ref{12-1} (resp. (3) in Theorem \ref{12}) can be replaced with `` $\K_{\mathrm{tac}}(R\text{-}\p)$ and $\K_{\mathrm{tac}}(R\text{-}\inj)$ are triangle equivalent'', where $\K_{\mathrm{tac}}(R\text{-}\p)$ and $\K_{\mathrm{tac}}(R\text{-}\inj)$ are as above.}
\end{rem}

\bigskip \centerline {\bf ACKNOWLEDGEMENTS}
We thank Lars Winther Christensen, Li Liang and Peder Thompson for conversations and comments on an early draft of this
paper. The work is partially supported by the National Natural Science Foundation
of China (Grant no. 12361008, 12061061) and the Foundation for
Innovative Fundamental Research Group Project of Gansu Province (Grant no. 23JRRA684),
as well as the grant PID 2020-113206GB-100 funded by  MCIN/AEI 10.13039/501100011033 as well as the grant 22004/PI/22 founded by the Fundaci\'{o}n S\'{e}neca.

\end{document}